\documentclass[11pt]{article}
\usepackage[margin=1in]{geometry}
\usepackage{amsfonts,amsmath,amssymb}
\usepackage{algpseudocode}

\PassOptionsToPackage{obeyspaces}{url}
\usepackage[colorlinks=true,citecolor=blue,urlcolor=blue,linkcolor=blue,bookmarksopen=true]{hyperref}
\usepackage{amsthm}
\usepackage{tikz}

\usepackage{etoolbox}
\patchcmd{\thebibliography}{\leftmargin\labelwidth}{\leftmargin\labelwidth\addtolength\itemsep{-0.1\baselineskip}}{}{}

\author{Christopher Cox\thanks{Department of Mathematics, Iowa State University, Ames, IA, USA. \texttt{cocox@iastate.edu}. Supported in part through NSF RTG Grant DMS-1839918.}
\and Ryan R.\ Martin\thanks{Department of Mathematics, Iowa State University, Ames, IA, USA. \texttt{rymartin@iastate.edu}. Supported in part through Simons Collaboration Grants \#353292 and \#709641.}}

\title{Counting paths, cycles and blow-ups in planar graphs}
\date{}

\usepackage[nameinlink,sort]{cleveref}

\newtheorem{theorem}{Theorem}
\newtheorem{lemma}[theorem]{Lemma}
\newtheorem{corollary}[theorem]{Corollary}

\newtheorem{conj}[theorem]{Conjecture}
\crefname{conj}{conjecture}{conjectures}

\newtheorem{claim}[theorem]{Claim}
\crefname{claim}{claim}{claims}

\newtheorem{prop}[theorem]{Proposition}
\crefname{prop}{proposition}{propositions}

\theoremstyle{definition}

\newtheorem{defn}[theorem]{Definition}
\crefname{defn}{definition}{definitions}

\newtheorem{remark}[theorem]{Remark}
\crefname{remark}{remark}{remarks}

\newtheorem{question}[theorem]{Question}
\crefname{question}{question}{questions}

\crefname{enumi}{part}{parts}

\numberwithin{theorem}{section}

\makeatletter
\DeclareRobustCommand{\crefnosort}[1]{%
    \begingroup\@cref@sortfalse\cref{#1}\endgroup
}
\DeclareRobustCommand{\Crefnosort}[1]{%
    \begingroup\@cref@sortfalse\Cref{#1}\endgroup
}
\makeatother

\newcommand*{\eqdef}{\stackrel{\mbox{\normalfont\tiny{def}}}{=}}        
\newcommand*{\abs}[1]{\lvert #1\rvert}                
\newcommand*{\abss}[1]{\bigl\lvert #1\bigr\rvert}     
\newcommand*{\absss}[1]{\biggl\lvert #1\biggr\rvert}  

\renewcommand*{\epsilon}{\varepsilon}       
\newcommand*{\R}{\mathbb{R}}                
\newcommand*{\Z}{\mathbb{Z}}                
\newcommand*{\wtilde}[1]{\widetilde{#1}}    
\newcommand*{\mcal}[1]{\mathcal{#1}}        
\newcommand*{\mbf}[1]{\mathbf{#1}}          

\newcommand*{\tp}[1]{\textproc{#1}}
\newcommand*{\plan}{\mathcal{P}}
\newcommand*{\optb}{\beta}
\newcommand*{\optp}{\rho}
\newcommand*{\sg}[1]{G_{#1}}
\newcommand*{\gcl}[1]{\mathcal{G}_{#1}}
\newcommand*{\bup}[2]{{#1\{#2\}}}
\DeclareMathOperator{\Opt}{Opt}
\DeclareMathOperator{\numb}{\mathbf{N}}
\DeclareMathOperator{\cp}{\mathbf{C}}

\let\Pr\relax\DeclareMathOperator*{\Pr}{\mathbf{Pr}}    
\DeclareMathOperator{\supp}{supp}   
\DeclareMathOperator{\Aut}{Aut}     

\begin{document}
\maketitle

\begin{abstract}
    For a planar graph $H$, let $\numb_\plan(n,H)$ denote the maximum number of copies of $H$ in an $n$-vertex planar graph.
    In this paper, we prove that $\numb_\plan(n,P_7)\sim{4\over 27}n^4$, $\numb_\plan(n,C_6)\sim(n/3)^3$, $\numb_\plan(n,C_8)\sim(n/4)^4$ and $\numb_\plan(n,\bup{K_4}{1})\sim(n/6)^6$, where $\bup{K_4}{1}$ is the $1$-subdivision of $K_4$.
    In addition, we obtain significantly improved upper bounds on $\numb_\plan(n,P_{2m+1})$ and $\numb_\plan(n,C_{2m})$ for $m\geq 4$.
    For a wide class of graphs $H$, the key technique developed in this paper allows us to bound $\numb_\plan(n,H)$ in terms of an optimization problem over weighted graphs.
\end{abstract}

\section{Introduction}

In this paper, we use standard graph theory definitions and notation (c.f.~\cite{west_graph}): $P_n$, $C_n$ and $K_n$ denote the path, cycle and clique on $n$ vertices, respectively.
The complete bipartite graph with parts of size $a$ and $b$ is denoted by $K_{a,b}$.
We use also standard big-oh and little-oh notation.

For graphs $G$ and $H$, let $\numb(G,H)$ denote the number of (unlabeled) copies of $H$ in $G$.
For a collection of graphs $\mcal G$ and a positive integer $n$, define
\[
    \numb_{\mcal G}(n,H)\eqdef\max\bigl\{\numb(G,H):G\in\mcal G,\ \abs{V(G)}=n\bigr\}.
\]
In this paper, we are concerned with asymptotically determining $\numb_\plan(n,H)$ for various graphs $H$, where $\plan$ is the set of all planar graphs.
\medskip

The study of $\numb_\plan(n,H)$ was initiated by Hakimi and Schmeichel~\cite{hakimi_cycles}, who determined $\numb_\plan(n,C_3)$ and $\numb_\plan(n,C_4)$ exactly.
Alon and Caro~\cite{alon_biclique} continued this study by determining $\numb_\plan(n,K_{2,k})$ exactly for all $k$; in particular, they determined $\numb_\plan(n,P_3)$.
Gy\H{o}ri et al.~\cite{gyori_p4} later gave the exact value for $\numb_\plan(n,P_4)$, and the same authors determined $\numb_\plan(n,C_5)$ in~\cite{gyori_c5}.
Generalizations of some of these results to other surfaces were established by Huynh, Joret and Wood~\cite{huynh_surface}.

The main driving force behind this manuscript is a recent conjecture of Ghosh et al.~\cite{ghosh_planarp5} which posits that
\begin{equation}\label{eqn:ghosh}
    \numb_\plan(n,P_{2m+1})=4m\biggl({n\over m}\biggr)^{m+1}+O(n^m)\qquad\text{for all }m\geq 2;
\end{equation}
the authors construct graphs which meet the lower bound for all $m\geq 2$, and they prove the case of $m=2$, showing that $\numb_\plan(n,P_5)=n^3+O(n^2)$.
We make steps toward this conjecture by proving:
\begin{theorem}\label{thm:oddpath}
    The following hold:
    \begin{alignat*}{2}
        \numb_\plan(n,P_7) &={4\over 27}\cdot n^{4}+O(n^{4-1/5}), &&\text{and}\\
        \numb_\plan(n,P_{2m+1}) &\leq {n^{m+1}\over2\cdot(m-1)!}+O(n^{m+4/5})\quad &&\text{for all }m\geq 4.
    \end{alignat*}
\end{theorem}
This, in particular, establishes the $m=3$ case of Ghosh et al.'s conjecture, albeit with a worse error-term than predicted.
Prior to this result, the best general upper bound that we are aware of is
\[
    \numb_\plan(n,P_{2m+1})\leq{(6n)^{m+1}\over 2}\qquad\text{for all }m\geq 3,
\]
though this bound does not appear to be in the literature.
\medskip

The methods used to prove this result extend to even cycles.

\begin{theorem}\label{thm:evencycle}
    The following hold:
    \begin{alignat*}{2}
        \numb_\plan(n,C_6) &=\biggl({n\over 3}\biggr)^3+O(n^{3-1/5}),\\
        \numb_\plan(n,C_8) &=\biggl({n\over 4}\biggr)^4+O(n^{4-1/5}),\quad &&\text{and}\\
        \numb_\plan(n,C_{2m}) &\leq {n^m\over m!}+O(n^{m-1/5}) &&\text{for all }m\geq 5.
    \end{alignat*}
\end{theorem}
Prior to this result, the best general upper bound that we are aware of is
\[
    \numb_\plan(n,C_{2m})\leq{(6n)^m\over 4m}\qquad\text{for all }m\geq 3.
\]

We present also new proofs of some known results.
\begin{theorem}\label{thm:alreadyknown}
    The following hold:
    \begin{enumerate}
        \item\label{known:p5} $\displaystyle\numb_\plan(n,P_5)=n^3+O(n^{14/5})$,\hfill(Ghosh et al.~\cite{ghosh_planarp5})
        \item\label{known:c4} $\displaystyle\numb_\plan(n,C_4)={n^2\over 2}+O(n^{9/5})$,\hfill(Hakimi--Schmeichel~\cite{hakimi_cycles})
        \item\label{known:k2k} $\displaystyle\numb_\plan(n,K_{2,k})={n^k\over k!}+O(n^{k-1+16/(k+8)})$\qquad for $k\geq 9$.\hfill(Alon--Caro~\cite{alon_biclique})
    \end{enumerate}
\end{theorem}
Although these results are already known and our error-terms are worse than those attained in the original papers, these results demonstrate the strength of the method developed in this paper.
Indeed, after applying one of a trio of general reduction lemmas (discussed in \Cref{sec:reduce}), each of these results follow in about one to two lines.
Furthermore, our results actually apply to a wider class of graphs than just planar graphs, namely the class of graphs which have linearly many edges and have no copy of $K_{3,3}$.
\medskip

Beyond odd paths and even cycles, our methods allow us to tackle particular blow-ups of graphs.

\begin{defn}\label{defn:blowup}
    Let $H=(V,E)$ be a graph and let $k$ be a positive integer.
    The \emph{$k$-edge-blow-up of $H$} is the graph $\bup{H}{k}$, which is formed by replacing every edge $xy\in E$ by an independent set of size $k$ and connecting each of these $k$ new vertices to both $x$ and $y$.
\end{defn}
For example, $\bup{C_m}{1}=C_{2m}$ for $m\geq 3$ and $\bup{K_2}{k}=K_{2,k}$ for $k\geq 1$.
We note that the graph $\bup{C_m}{\ell}$ where $\ell=\lfloor{n-m\over m}\rfloor$ realizes the lower-bound in \cref{eqn:ghosh}.
\medskip

Alon and Caro~\cite{alon_biclique} determined $\numb_\plan(n,\bup{K_2}{k})$ exactly for all $k$; we extend this to the other two planar cliques by showing:
\begin{theorem}\label{thm:planarclique}
    For all positive integers $k$,
    \begin{align*}
        \numb_\plan(n,\bup{K_3}{k}) &= {1\over(k!)^3}\biggl({n\over 3}\biggr)^{3k}+O(n^{3k-k/(k+4)}),\quad\text{and}\\
        \numb_\plan(n,\bup{K_4}{k}) &= {1\over(k!)^6}\biggl({n\over 6}\biggr)^{6k}+O(n^{6k-k/(k+4)}).
    \end{align*}
\end{theorem}

In general, it is not difficult to show that $\numb_\plan(n,\bup{H}{k})=\Theta(n^{km})$ if $H$ is a planar graph on $m$ edges and $k\cdot\delta(H)\geq 2$.
Indeed, the graph $\bup{H}{\ell}$ where $\ell=\lfloor{n-\abs{V(H)}\over m}\rfloor$ shows that
\begin{equation}\label{eqn:blowuplower}
    \numb_\plan(n,\bup{H}{k})\geq{\ell\choose k}^m={1\over (k!)^m}\biggl({n\over m}\biggr)^{km}-O(n^{km-1}),
\end{equation}
and it is an exercise to bound
\begin{equation}\label{eqn:blowupper}
    \numb_\plan(n,\bup{H}{k})\leq{(6n)^{km}\over \abs{\Aut H}\cdot (k!)^m},
\end{equation}
where $\Aut H$ is the automorphism group of $H$.
The key step in the proof of this upper bound is the content of \Cref{prop:easyupper}.
In this paper, we significantly improve the leading constant in the upper-bound.

\begin{theorem}\label{thm:blowupbound}
    Let $H$ be a planar graph on $m$ edges and let $k$ be a positive integer.
    If either
    \begin{itemize}
        \item $k\cdot\bigl(\delta(H)-1\bigr)\geq 2$, or
        \item $\delta(H)=1$ and $k\geq 9$,
    \end{itemize}
    then
    \[
        \numb_\plan(n,\bup{H}{k})\leq{n^{km}\over (km)!}+o(n^{km}).
    \]
\end{theorem}
Compare this result to the na\"ive bounds in \cref{eqn:blowuplower,eqn:blowupper}.
In fact, provided $k$ is sufficiently large, we are able to asymptotically pin down $\numb_\plan(n,\bup{H}{k})$.
\begin{theorem}\label{thm:blowupexact}
    Let $H$ be a planar graph on $m$ edges and let $k$ be a positive integer.
    If either
    \begin{itemize}
        \item $\delta(H)\geq 2$ and $k\geq{\log(m+1)\over m\log(1+1/m)}$, or
        \item $\delta(H)=1$ and $k\geq\max\bigl\{9,{\log(m+1)\over m\log(1+1/m)}\bigr\}$,
    \end{itemize}
    then
    \[
        \numb_\plan(n,\bup{H}{k})={1\over(k!)^m}\biggl({n\over m}\biggr)^{km}+o(n^{km}).
    \]
\end{theorem}
The requirement that $k\geq{\log(m+1)\over m\log(1+1/m)}$ in the above theorem is necessary for some graphs $H$.
As an example, let $I$ denote the skeleton of the icosahedron and let $I^-$ denote the graph formed by deleting any edge from $I$.
Since $\abs{E(I^-)}=29$ and $\delta(I^-)=4$, \Cref{thm:blowupexact} implies that $\numb_\plan(n,\bup{I^-}{k})\sim{1\over (k!)^{29}}\bigl({n\over 29}\bigr)^{29k}$ for all $k\geq 4$.
However, for $k\in\{1,2,3\}$, the graph $\bup{I}{\ell}$ where $\ell=\lfloor{n-12\over 30}\rfloor$ realizes
\[
    \numb_\plan(n,\bup{I^-}{k})\geq 30{\ell\choose k}^{29}\sim{30\over(k!)^{29}}\biggl({n\over 30}\biggr)^{29k}>{1.57\over(k!)^{29}}\biggl({n\over 29}\biggr)^{29k},
\]
since $\numb(I,I^-)=30$.
The icosahedron is not unique in this regard (see \Cref{prop:edgetrans}).
\medskip

The paper is organized as follows.
In \Cref{sec:reduce}, we present the key contribution of this paper: a trio of reduction lemmas from which all of our results follow.
\Cref{sec:reduceproof} contains the proofs of these reduction lemmas.
We then, in \Cref{sec:paths}, use these reduction lemmas to prove \Cref{thm:oddpath} and \cref{known:p5} of \Cref{thm:alreadyknown}.
In \Cref{sec:blowup}, we establish \Cref{thm:evencycle,thm:planarclique,thm:blowupbound,thm:blowupexact} along with \cref{known:c4,known:k2k} of \Cref{thm:alreadyknown}.
We conclude with a list of remarks and conjectures in \Cref{sec:remarks}.

\subsection{Notation and preliminaries}
We use standard graph theory definitions and notation (c.f.~\cite{west_graph}).
For a graph $G$, we use $V(G)$ and $E(G)$ to denote its vertex-set and edge-set, respectively.
When the graph is understood, we omit the parenthetical and simply write $V$ and $E$.

For $v\in V(G)$, we write $N_G(v)$ to denote the neighborhood of $v$ in $G$ and $\deg_G(v)\eqdef\abs{N_G(v)}$ to denote the degree of $v$ in $G$.
For vertices $u,v\in V(G)$, we write $\deg_G(u,v)\eqdef\abs{N_G(u)\cap N_G(v)}$ to denote the co-degree of $u$ and $v$ in $G$.
When the graph $G$ is understood, we omit the subscript.

For positive integers $m\leq n$, we write $[n]$ to denote the set $\{1,\dots,n\}$ and write $[m,n]$ to denote the set $\{m,\dots,n\}$.
For a set $X$, we use $(X)_n$ to denote the set of tuples $(x_1,\dots,x_n)\in X^n$ with $x_1,\dots,x_n$ distinct; this notation mirrors that of the falling-factorial.
Finally, we will often write $xy$ to denote the set $\{x,y\}$ for notational convenience.
\medskip

We require a special case of the Karush--Kuhn--Tucker (KKT) conditions (c.f.~\cite[Corollaries 9.6 and 9.10]{guler_opt}) in order to prove \Cref{lem:aequalb,lem:regularity}.
\begin{theorem}[Special case of the KKT conditions]\label{KKT}
    Let $f\colon\R^n\to\R$ be a continuously differentiable function and consider the optimization problem
    \[
        \begin{array}{cl}
            \max & f(\mbf x)\\
            \text{s.t.} & \sum_ix_i=1\\
                        & x_1,\dots,x_n\geq 0.
        \end{array}
    \]
    If $\mbf x^*$ achieves this maximum, then there is some $\lambda\in\R$ such that, for each $i\in[n]$, either
    \[
        x_i^*=0,\qquad\text{or}\qquad{\partial f\over\partial x_i}(\mbf x^*)=\lambda.
    \]
\end{theorem}

\section{The key reduction lemmas}\label{sec:reduce}

Aside from the bounds stated in the introduction, the main contribution of this paper is the technique used in their proofs.

For graphs $G,H$, let $\cp(G,H)$ denote the set of (unlabeled) copies of $H$ in $G$; so $\abs{\cp(G,H)}=\numb(G,H)$.
For a finite set $X$, we abbreviate $\cp(X,H)\eqdef\cp(K_X,H)$, where $K_X$ is the clique on vertex-set $X$; in other words, $\cp(X,H)$ is the set of all copies $H'$ of $H$ with $V(H')\subseteq X$.

The following definition lays out the key quantities used throughout this paper.

\begin{defn}\label{defn:opt}
    Fix a finite set $X$ and let $\mu$ be a probability mass on ${X\choose 2}$.
    We define the following quantities:
    \begin{enumerate}
        \item For $x\in X$, define
            \[
                \bar\mu(x)\eqdef\sum_{y\in X\setminus\{x\}}\mu(xy),
            \]
            which is the probability that an edge sampled from $\mu$ is incident to $x$.
            It can also be thought of as the weighted degree of $x$.
            Note that $\sum_{x\in X}\bar\mu(x)=2$ thanks to the handshaking lemma.
        \item For an integer $m\geq 2$, define
            \begin{align*}
                \optp(\mu;m) &\eqdef\sum_{\mbf x\in(X)_m}\bar\mu(x_1)\biggl(\prod_{i=1}^{m-1}\mu(x_ix_{i+1})\biggr)\bar\mu(x_m),\quad\text{and}\\
                \optp(m) &\eqdef\sup\biggl\{\optp(\mu;m):\supp\mu\subseteq{X\choose 2}\text{ for some finite set }X\biggr\}.
            \end{align*}
            The quantity $\optp(\mu;m)$ is essentially the probability that, upon independently sampling edges $e_1,\dots,e_{m+1}$ from $\mu$, the edges $e_2,\dots,e_m$ form a copy of $P_m$, $e_1$ is incident to the first vertex of this path and $e_{m+1}$ is incident to the last vertex of this path (see the proof of \Cref{thm:probp} for a more precise interpretation).
        \item For a subgraph $G\subseteq K_X$, define
            \[
                \mu(G)\eqdef\prod_{e\in E(G)}\mu(e),
            \]
            which is essentially the probability that $\abs{E(G)}$ edges sampled independently from $\mu$ form the edge-set of $G$.
        \item For a graph $H$ with no isolated vertices and a positive integer $k$, define
            \begin{align*}
                \optb(\mu;H,k) &\eqdef\sum_{H'\in \cp(X,H)}\mu(H')^k,\quad\text{and}\\
                \optb(H,k) &\eqdef\sup\biggl\{\optb(\mu;H,k):\supp\mu\subseteq{X\choose 2}\text{ for some finite set }X\biggr\}.
            \end{align*}
            The quantity $\optb(\mu;H,k)$ is essentially the probability that $k\cdot\abs{E(H)}$ edges sampled independently from $\mu$ form a copy of $H$ wherein each edge has multiplicity $k$ (see the proof of \Cref{thm:probb} for a more precise interpretation).
    \end{enumerate}
\end{defn}

While we are primarily concerned with planar graphs, our results apply to a much broader class of graphs.
\begin{defn}
    For any fixed $C>0$, the collection of graphs $\gcl C$ is defined as follows:
    $G\in\gcl C$ if and only if
    \begin{enumerate}
        \item $G$ has no copy of $K_{3,3}$, and
        \item Every subgraph $H\subseteq G$ satisfies $\abs{E(H)}\leq C\cdot\abs{V(H)}$.
    \end{enumerate}
\end{defn}

Observe that $\gcl{C_1}\subseteq\gcl{C_2}$ if $C_1\leq C_2$ and that $\plan\subseteq\gcl 3$.
Furthermore, observe that $\bup{H}{k}\in\gcl 2$ for any graph $H$ and any positive integer $k$.
In each of the results discussed in the introduction, $\plan$ can be replaced by $\gcl C$ for any $C\geq 2$ (due to monotonicity, all of our upper-bounds hold for any $C>0$, but the lower-bound constructions require $C\geq 2$).
\medskip

We quickly remark that our results apply to an even wider class of graphs than $\gcl C$, though we avoid this more general situation for the sake of readability.
We discuss these further generalizations in \Cref{sec:remarks}.
\medskip

For paths of odd order, we show:
\begin{lemma}[Reduction lemma for odd paths]\label{lem:tumorp}
    If $m\geq 2$, then
    \[
        \numb_{\gcl C}(n,P_{2m+1})\leq{\optp(m)\over 2}\cdot n^{m+1}+O(n^{m+4/5}),
    \]
    where the implicit constant in the big-oh notation depends on $m$ and $C$.
\end{lemma}

For general edge-blow-ups, we prove:
\begin{lemma}[Reduction lemma for edge-blow-ups]\label{lem:tumorb}
    Let $H$ be a graph on $m$ edges and let $k$ be a positive integer.
    If $k\cdot\bigl(\delta(H)-1\bigr)\geq 2$, then
    \[
        \numb_{\gcl C}(n,\bup{H}{k})\leq{\optb(H,k)\over(k!)^m}\cdot n^{km}+O(n^{km-k/(k+4)}).
    \]
    If $\delta(H)=1$ and $k\geq 9$, then
    \[
        \numb_{\gcl C}(n,\bup{H}{k})\leq{\optb(H,k)\over(k!)^m}\cdot n^{km}+O(n^{km-1+16/(k+8)}).
    \]
    In both cases, the implicit constant in the big-oh notation depends on $H$, $k$ and $C$.
\end{lemma}

Recall that $C_{2m}=\bup{C_m}{1}$ for $m\geq 3$ and that $C_4=\bup{K_2}{2}$.
Unfortunately, since $\delta(C_m)=2$ for $m\geq 3$ and $\delta(K_2)=1$, we cannot apply \Cref{lem:tumorb} to these graphs.
However, with a slightly different approach, we can obtain exactly this extension.

\begin{lemma}[Reduction lemma for even cycles]\label{lem:tumorc}
    The following hold:
    \begin{alignat*}{2}
        \numb_{\gcl C}(n,C_4) &\leq{\optb(K_2,2)\over 2}\cdot n^2+O(n^{2-1/5}), && \text{and}\\
        \numb_{\gcl C}(n,C_{2m}) &\leq\optb(C_m,1)\cdot n^m+O(n^{m-1/5})\quad &&\text{for }m\geq 3,
    \end{alignat*}
    where the implicit constant in the big-oh notation depends on $m$ and $C$.
\end{lemma}

Note that there is still a gap in the reduction lemmas when it comes to $\bup{K_2}{k}=K_{2,k}$, which we can handle only if $k=2$ or $k\geq 9$; we suspect that this gap can be closed.
Granted, at least when dealing with $\plan$, this result is already superseded by the results of Alon and Caro~\cite{alon_biclique}.
However, we believe that the obvious reduction lemma holds for $\bup{H}{k}$ provided $k\cdot\delta(H)\geq 2$, though we do not currently have a proof.
\medskip

While the individual details in each of these reduction lemmas differ, the underlying philosophy is the same.
The key idea is to show that the vast majority of the copies of $H$ in $G$ interact predictably with the largest degree vertices of $G$.
This being the case, we then argue that $G$ can be suitably approximated by an edge-blow-up of some graph, possibly where each edge is blown up by different amounts.
The probability masses $\mu$ discussed in \Cref{defn:opt} are a compact way to represent these edge-blow-ups which approximate $G$.

\subsection{Proofs of the reduction lemmas}\label{sec:reduceproof}
In this section, we prove \Cref{lem:tumorp,lem:tumorb,lem:tumorc}.
The approach to the lemmas is very similar, yet each requires separate analysis.
\medskip

We begin by presenting a simple proposition, pieces of which are used in each proof.

\begin{prop}\label{prop:codegreebound}
    Let $G=(V,E)\in\gcl C$ be a graph on $n$ vertices.
    For $\epsilon>0$, define $\wtilde V\eqdef\{v\in V:\deg(v)\geq\epsilon n\}$.
    Then,
    \[
        \abs{\wtilde V}\leq {2C\over\epsilon},\qquad\text{and}\qquad \sum_{uv\in{\wtilde V\choose 2}}\deg(u,v) \leq n+4\biggl({C\over\epsilon}\biggr)^4.
    \]
\end{prop}
\begin{proof}
    We begin by observing that
    \[
        \epsilon n\cdot\abs{\wtilde V}\leq\sum_{v\in\wtilde V}\deg(v)\leq\sum_{v\in V}\deg(v)=2\abs E\leq 2Cn\quad\implies\quad \abs{\wtilde V}\leq{2C\over\epsilon}.
    \]

    For notational convenience set $\wtilde E\eqdef{\wtilde V\choose 2}$ and $S\eqdef\sum_{uv\in\wtilde E}\deg(u,v)$.
    Since $G$ has no copy of $K_{3,3}$, we know that $\abss{N(u)\cap N(v)\cap N(w)}\leq 2$ for any distinct $u,v,w\in V$.
    Hence, we can apply the second Bonferroni inequality to bound
    \begin{align*}
        n &\geq \absss{\bigcup_{uv\in\wtilde E}\bigl(N(u)\cap N(v)\bigr)}\geq\sum_{uv\in\wtilde E}\abss{N(u)\cap N(v)}-\sum_{\{uv,wz\}\in{\wtilde E\choose 2}}\abss{N(u)\cap N(v)\cap N(w)\cap N(z)}\\
          &\geq S-2{\abs{\wtilde E}\choose 2}\geq S-{1\over 4}\abs{\wtilde V}^4\geq S-4\biggl({C\over\epsilon}\biggr)^4,
    \end{align*}
    which proves the proposition.
\end{proof}

\paragraph{Reduction lemma for odd paths.}

\begin{proof}[Proof of \Cref{lem:tumorp}]
    Fix $G=(V,E)\in\gcl C$ on $n$ vertices and fix $\mbf v=(v_1,\dots,v_m)\in(V)_m$.
    Label $V(P_{2m+1})=\{p_1,p_2,\dots,p_{2m+1}\}$ in consecutive order and consider the copies of $P_{2m+1}$ in $G$ wherein $v_i$ plays the role of vertex $p_{2i}$.
    Observe that there are then at most $\deg(v_1)$ choices for the image of $p_1$, at most $\deg(v_i,v_{i+1})$ choices for the image of $p_{2i+1}$ for $i\in[m-1]$ and at most $\deg(v_m)$ choices for the image of $p_{2m+1}$.
    Thus, there are at most
    \[
        D(\mbf v)\eqdef\deg(v_1)\biggl(\prod_{i=1}^{m-1}\deg(v_i,v_{i+1})\biggr)\deg(v_m)
    \]
    copies of $P_{2m+1}$ in $G$ wherein $v_i$ plays the role of vertex $p_{2i}$ and so we can bound
    \[
        \numb(G,P_{2m+1})\leq{1\over 2}\sum_{\mbf v\in(V)_m}D(\mbf v).
    \]

    Fix $\epsilon=\epsilon(n)>0$ to be chosen later and define the set
    \[
        \wtilde E\eqdef\biggl\{uv\in{V\choose 2}:\deg(u,v)\geq\epsilon n\biggr\}.
    \]
    The set $\wtilde E$ induces a graph $\wtilde G$ with vertex-set $\wtilde V\subseteq V$.
    Certainly if $v\in\wtilde V$, then $\deg(v)\geq\epsilon n$ and so $\abs{\wtilde V}\leq 2C/\epsilon$, thanks to \Cref{prop:codegreebound}.

    Next define
    \begin{align*}
        \wtilde P_m &\eqdef\bigl\{\mbf v\in (V)_m: v_iv_{i+1}\in\wtilde E\text{ for all $i\in[m-1]$}\bigr\},\quad\text{and}\\
        \wtilde M &\eqdef\sum_{\mbf v\in\wtilde P_m}D(\mbf v).
    \end{align*}
    We aim to show that $\numb(G,P_{2m+1})\approx\wtilde M/2$.

    For any $u,v\in V$, we have $\deg(u,v)\leq\min\{\deg(u),\deg(v)\}$, so for any $\mbf v\in(V)_m$, we can bound
    \begin{align*}
        D(\mbf v) &\leq \biggl(\prod_{i=1}^j\deg(v_i)\biggr)\cdot\deg(v_j,v_{j+1})\cdot\biggl(\prod_{i=j+1}^m\deg(v_i)\biggr)\quad\text{for all $j\in[m-1]$}\\
        \implies D(\mbf v) &\leq\biggl(\min_{i\in[k-1]}\deg(v_i,v_{i+1})\biggr)\prod_{i=1}^m\deg(v_i).
    \end{align*}
    We can therefore bound
    \begin{align*}
        2\cdot\numb(G,P_{2m+1})-\wtilde M &\leq\sum_{\mbf v\in (V)_m\setminus\wtilde P_m}D(\mbf v) \leq\sum_{\mbf v\in(V)_m\setminus\wtilde P_m}\biggl(\min_{i\in[k-1]}\deg(v_i,v_{i+1})\biggr)\prod_{i=1}^m\deg(v_i)\\
                                          &\leq\sum_{\mbf v\in (V)_m\setminus\wtilde P_m}\epsilon n\cdot \prod_{i=1}^m\deg(v_i)\leq\epsilon n\cdot\sum_{v_1,\dots,v_m\in V}\ \prod_{i=1}^m\deg(v_i)\\
                                          &=\epsilon n\cdot\biggl(\sum_{v\in V}\deg(v)\biggr)^m\leq \epsilon n\cdot (2Cn)^m=O\bigl(\epsilon n^{m+1}\bigr).
    \end{align*}

    Set $U=\{v\in V\setminus\wtilde V:\abs{N(v)\cap \wtilde V}\geq 3\}$ and define the subgraph $G'=(V',E')$ of $G$ as follows:
    \begin{itemize}
        \item Delete all vertices in $U$, and
        \item Delete all vertices $v\in V\setminus\wtilde V$ for which $N(v)\cap\wtilde V=\varnothing$, and
        \item Delete all edges induced by $\wtilde V$, and
        \item Delete all edges induced by $V\setminus\wtilde V$.
    \end{itemize}
    Since $G$ has no copy of $K_{3,3}$,
    \[
        \abs U\leq 2{\abs{\wtilde V}\choose 3}\leq 2{2C/\epsilon\choose 3}=O(\epsilon^{-3}).
    \]

    For $\mbf v\in (\wtilde V)_m$, define
    \begin{align*}
        D'(\mbf v) &\eqdef\deg_{G'}(v_1)\biggl(\prod_{i=1}^{m-1}\deg_{G'}(v_i,v_{i+1})\biggr)\deg_{G'}(v_m),\quad\text{and}\\
        \wtilde M' &\eqdef\sum_{\mbf v\in\wtilde P_m}D'(\mbf v).
    \end{align*}
    For $v\in\wtilde V$, observe that
    \[
        {\deg_{G'}(v)\over\deg_G(v)}=1-{\deg_G(v)-\deg_{G'}(v)\over\deg_G(v)}\geq 1-{O(\epsilon^{-1})+O(\epsilon^{-3})\over\epsilon n}=1-O\biggl({1\over\epsilon^4 n}\biggr).
    \]
    Similarly, for $uv\in\wtilde E$,
    \[
        {\deg_{G'}(u,v)\over\deg_G(u,v)}\geq 1-{O(\epsilon^{-1})+O(\epsilon^{-3})\over \epsilon n}=1-O\biggl({1\over\epsilon^4 n}\biggr).
    \]
    Therefore, for any $\mbf v\in\wtilde P_m$, we have
    \[
        {D'(\mbf v)\over D(\mbf v)} = {\deg_{G'}(v_1) \over\deg_G(v_1)}\biggl(\prod_{i=1}^{m-1}{\deg_{G'}(v_i,v_{i+1})\over\deg_G(v_i,v_{i+1})}\biggr){\deg_{G'}(v_m)\over\deg_G(v_m)}\geq 1-O\biggl({1\over\epsilon^4 n}\biggr),
    \]
    and so
    \[
        \wtilde M' \geq \biggl(1-O\biggl({1\over\epsilon^4 n}\biggr)\biggr)\wtilde M.
    \]

    Next, we can partition $V'\setminus\wtilde V=U_1\cup U_2$ where $U_i=\{v\in V'\setminus\wtilde V:\abs{N_{G'}(v)}=i\}$.
    We claim that we may suppose $U_1=\varnothing$.
    Indeed, suppose that $x\in U_1$ and that $xu\in E(G')$, so $u\in\wtilde V$.
    Consider selecting any $v$ such that $uv\in\wtilde E$ and introducing the edge $xv$.
    (Note that if $G$ was planar to begin with, then $G'$ is still planar after this modification.)
    Observe that $\wtilde M'$ can only increase under this operation and so we may suppose that $U_1=\varnothing$.

    Thus, set $S\eqdef\sum_{uv\in{\wtilde V\choose 2}}\deg_{G'}(u,v)$ and let $\mu$ be the probability mass on ${\wtilde V\choose 2}$ defined by $\mu(uv)=\deg_{G'}(u,v)/S$.
    Since $V'=\wtilde V\cup U_2$, and $G'$ has no edges induced by $\wtilde V$, we observe that $S=\abs{U_2}\leq n$.
    Furthermore, for any $v\in\wtilde V$, we have $\bar\mu(v)=\deg_{G'}(v)/S$.
    Therefore,
    \begin{align*}
        \wtilde M' &= \sum_{\mbf v\in\wtilde P_m}D'(\mbf v)\leq\sum_{\mbf v\in(\wtilde V)_m}D'(\mbf v)= S^{m+1}\cdot\sum_{\mbf v\in(\wtilde V)_m}\bar\mu(v_1)\biggl(\prod_{i=1}^{m-1}\mu(v_iv_{i+1})\biggr)\bar\mu(v_m)\\
                   &= \optp(\mu;m)\cdot S^{m+1}\leq \optp(\mu;m)\cdot n^{m+1}\leq\optp(m)\cdot n^{m+1}.
    \end{align*}
    Finally, selecting $\epsilon=n^{-1/5}$ yields
    \begin{align*}
        \numb(G,P_{2m+1}) &\leq {1\over 2}\wtilde M+O(\epsilon n^{m+1})\leq{1\over 2}\biggl(1+O\biggl({1\over\epsilon^4n}\biggr)\biggr)\optp(m)\cdot n^{m+1}+O(\epsilon n^{m+1})\\
                        &={\optp(m)\over 2}\cdot n^{m+1}+O(n^{m+4/5}).\qedhere
    \end{align*}
\end{proof}

Before moving on, we make a few remarks.
\begin{remark}
    It is not difficult to argue that for $m\geq 2$ and $C\geq 2$,
    \[
        \numb_{\gcl C}(n,P_{2m+1})\geq{\optp(m)\over 2}\cdot n^{m+1}-o(n^{m+1}),
    \]
    so \Cref{lem:tumorp} is asymptotically tight.
    Indeed, fix a finite set $X$ and a probability mass $\mu$ on ${X\choose 2}$.
    For a sufficiently large integer $n$, let $G$ be the edge-blow-up of $K_X$ formed by blowing up each edge $e\in{X\choose 2}$ into a set of size $\lfloor n\cdot\mu(e)\rfloor$.
    Then, one can show that $G\in\gcl 2$ and
    \[
        \numb(G,P_{2m+1})\geq{\optp(\mu;m)\over 2}\cdot n^{m+1}-O(n^m).
    \]
\end{remark}
\begin{remark}\label{rem:planarpath}
    For a finite set $X$ and a probability mass $\mu$ on ${X\choose 2}$, let $\sg\mu$ be the graph with vertex-set $X$ and edge-set $\supp\mu$.
    In the proof of \Cref{lem:tumorp}, if $G\in\plan$, then we can guarantee also that $G'\in\plan$, even after the modification that ensures $U_1=\varnothing$.
    Therefore, we can actually establish
    \[
        \numb_\plan(n,P_{2m+1})={\optp_\plan(m)\over 2}\cdot n^{m+1}+o(n^{m+1})\qquad\text{for }m\geq 2,
    \]
    where
    \[
        \optp_\plan(m)=\sup\bigl\{\optp(\mu;m):\sg\mu\in\plan\bigr\}.
    \]
    Although this refinement exists, we do not believe it to be helpful here.
    That is to say, we believe that $\optp_\plan(m)=\optp(m)$ for all $m\geq 2$.
\end{remark}

\paragraph{Reduction lemma for edge-blow-ups.}
We will need another simple proposition in order to establish \Cref{lem:tumorb}.
\begin{prop}\label{prop:easyupper}
    Let $H$ be a graph on $m$ edges and let $k$ be a positive integer.
    If $G=(V,E)$ is any graph and $k\cdot\delta(H)\geq 2$, then
    \[
        \sum_{H'\in\cp(V,H)}\ \prod_{xy\in E(H')}\deg_G(x,y)^k\leq{\bigl(2\abs E\bigr)^{km}\over\abs{\Aut H}}.
    \]
\end{prop}
\begin{proof}
    Since $k\cdot\delta(H)\geq 2$, we know that for any $x\in\R^+$ and $v\in V(H)$, we have
    \[
        1+x^{k\cdot\deg_H(v)/2}\leq (1+x)^{k\cdot\deg_H(v)/2}.
    \]
    Additionally, for $u\neq v\in V(G)$, observe that $\deg(u,v)\leq\min\{\deg(u),\deg(v)\}\leq\sqrt{\deg(u)\deg(v)}$.

    Using these two facts and translating between labeled and unlabeled copies of $H$, we can bound
    \begin{align*}
        \sum_{H'\in\cp(V,H)}\ \prod_{xy\in E(H')}\deg_G(x,y)^k &\leq \sum_{H'\in\cp(V,H)}\ \prod_{xy\in E(H')}\bigl(\deg_G(x)\deg_G(y)\bigr)^{k/2}\\
                                                               &= \sum_{H'\in\cp(V,H)}\ \prod_{x\in V(H')}\deg_G(x)^{k\cdot\deg_{H'}(x)/2}\\
                                                               &={1\over\abs{\Aut H}}\sum_{\substack{g\colon V(H)\to V\\ g\text{ injection}}}\ \prod_{v\in V(H)}\deg_G\bigl(g(v)\bigr)^{k\cdot\deg_H(v)/2}\\
                                                               &\leq{1\over\abs{\Aut H}}\sum_{g\colon V(H)\to V}\ \prod_{v\in V(H)}\deg_G\bigl(g(v)\bigr)^{k\cdot\deg_H(v)/2}.
    \end{align*}
    From here, we use the fact that
    \[
        \sum_{x_1,\dots,x_n\in X}\ \prod_{i=1}^n f_i(x_i)=\prod_{i=1}^n\biggl(\sum_{x\in X}f_i(x)\biggr)
    \]
    for any finite set $X$ and any functions $f_1,\dots,f_n\colon X\to\R$ in order to bound
    \begin{align*}
        \sum_{H'\in\cp(V,H)}\ \prod_{xy\in E(H')}\deg_G(x,y)^k &\leq{1\over\abs{\Aut H}}\prod_{v\in V(H)}\biggl(\sum_{x\in V}\deg_G(x)^{k\cdot\deg_H(v)/2}\biggr)\\
                                                               &\leq{1\over\abs{\Aut H}}\prod_{v\in V(H)}\biggl(\sum_{x\in V}\deg_G(x)\biggr)^{k\cdot\deg_H(v)/2} ={\bigl(2\abs E\bigr)^{km}\over\abs{\Aut H}}.\qedhere
    \end{align*}
\end{proof}

\begin{proof}[Proof of \Cref{lem:tumorb}]
    Fix $G=(V,E)\in\gcl C$ on $n$ vertices.
    Fix an injection $g\colon V(H)\to V$ and consider the of copies of $\bup{H}{k}$ in $G$ where, for each $v\in V(H)$, $g(v)$ plays the role of vertex $v$.
    For each $uv\in E(H)$, observe that there are at most ${\deg(g(u),g(v))\choose k}$ choices for the $k$ common neighbors of $u,v$ in $\bup{H}{k}$; thus there are at most
    \[
        \prod_{uv\in E(H)}{\deg\bigl(g(u),g(v)\bigr)\choose k}
    \]
    ways to extend $g$ to an embedding of $\bup{H}{k}$.
    In particular, we can bound
    \begin{align}
        \numb(G,\bup{H}{k}) &\leq\sum_{H'\in\cp(V,H)}\ \prod_{xy\in E(H')}{\deg(x,y)\choose k} \leq{1\over(k!)^m}\sum_{H'\in\cp(V,H)}\ \prod_{xy\in E(H')}\deg(x,y)^k.\label{eqn:codegupper}
    \end{align}
    \medskip

    Fix $\epsilon=\epsilon(n)>0$ to be chosen later and define
    \begin{align*}
        \wtilde V &\eqdef\{v\in V:\deg(v)\geq\epsilon n\},\quad\text{and}\\
        \wtilde M &\eqdef\sum_{H'\in\cp(\wtilde V,H)}\ \prod_{xy\in E(H')}\deg(x,y)^k.
    \end{align*}
    We claim that $\wtilde M$ approximates $\numb(G,\bup{H}{k})$.
    The proof of this fact depends heavily on the minimum degree of $H$, so we break the proof into two claims.
    \begin{claim}\label{claim:mindeg2}
        If $k\cdot\bigl(\delta(H)-1\bigr)\geq 2$, then
        \[
            \numb(G,\bup{H}{k})\leq {\wtilde M\over (k!)^m}+O(\epsilon^k n^{km}).
        \]
    \end{claim}
    \begin{proof}
        Set $\wtilde E={\wtilde V\choose 2}$.
        For $H'\in\cp(V,H)$, observe that $V(H')\subseteq\wtilde V$ if and only if $E(H')\subseteq\wtilde E$.
        Furthermore, observe that if $uv\notin\wtilde E$, then $\deg(u,v)<\epsilon n$.
        Using \cref{eqn:codegupper}, we begin by bounding
        \begin{align*}
            (k!)^m\cdot\numb(G,\bup{H}{k})-\wtilde M &\leq \sum_{\substack{H'\in\cp(V,H):\\ V(H')\not\subseteq\wtilde V}}\ \prod_{xy\in E(H')}\deg(x,y)^k=\sum_{\substack{H'\in\cp(V,H):\\ E(H')\not\subseteq\wtilde E}}\ \prod_{xy\in E(H')}\deg(x,y)^k\\
                                               &\leq \sum_{\substack{H'\in\cp(V,H):\\ E(H')\not\subseteq\wtilde E}}\ \sum_{e\in E(H')\setminus\wtilde E} (\epsilon n)^k\cdot\prod_{xy\in E(H')\setminus\{e\}}\deg(x,y)^k\\
                                               &\leq (\epsilon n)^k\cdot\sum_{e\in E(H)}\ \sum_{H''\in\cp(V,H-e)}\ \prod_{xy\in E(H'')}\deg(x,y)^k
        \end{align*}
        Now, for any $e\in E(H)$, we have $k\cdot\delta(H-e)\geq k\cdot\bigl(\delta(H)-1\bigr)\geq 2$, and so we can apply \Cref{prop:easyupper} to $H-e$ to bound
        \begin{align*}
            (k!)^m\cdot\numb(G,\bup{H}{k})-\wtilde M &\leq (\epsilon n)^k\cdot\sum_{e\in E(H)}{\bigl(2\abs E\bigr)^{k(m-1)}\over\abs{\Aut(H-e)}}\\
                                               &\leq \epsilon^k\cdot m\cdot(2C)^{k(m-1)}\cdot n^{km}=O(\epsilon^k n^{km}).\qedhere
        \end{align*}
    \end{proof}

    \begin{claim}\label{claim:mindeg1}
        If $\delta(H)=1$ and $k\geq 2$, then
        \[
            \numb(G,\bup{H}{k})\leq {\wtilde M\over (k!)^m}+O(\epsilon^{k/2} n^{km+1}).
        \]
    \end{claim}
    \begin{proof}
        The proof of this fact is very similar to the proof of \Cref{prop:easyupper}.
        For $u,v\in V$, certainly $\deg(u,v)\leq\min\{\deg(u),\deg(v)\}\leq\sqrt{\deg(u)\deg(v)}$.
        Thus, by applying \cref{eqn:codegupper}, we can bound
        \begin{align*}
            (k!)^m\numb(G,\bup{H}{k})-\wtilde M &= \sum_{\substack{H'\in\cp(V,H):\\ V(H')\not\subseteq\wtilde V}}\ \prod_{xy\in E(H')}\deg(x,y)^k\leq\sum_{\substack{H'\in\cp(V,H):\\ V(H)\not\subseteq\wtilde V}}\ \prod_{xy\in E(H')}\bigl(\deg(x)\deg(y)\bigr)^{k/2}\\
                                                     &=\sum_{\substack{H'\in\cp(V,H):\\ V(H')\not\subseteq\wtilde V}}\ \prod_{x\in V(H')}\deg(x)^{k\cdot\deg_{H'}(x)/2}\\
                                                     &\leq\sum_{\substack{H'\in\cp(V,H):\\ V(H')\not\subseteq\wtilde V}}\ \sum_{y\in V(H')\setminus\wtilde V}(\epsilon n)^{k\cdot\deg_{H'}(y)/2}\prod_{x\in V(H')\setminus\{y\}}\deg(x)^{k\cdot\deg_{H'}(x)/2}
        \end{align*}

        Next, by translating between labeled and unlabeled copies of $H$, we continue to bound
        \begin{align*}
            (k!)^m\numb(G,\bup{H}{k})-\wtilde M &\leq{1\over\abs{\Aut H}}\sum_{v\in V(H)}\ \sum_{\substack{g\colon V(H)\to V\\ g(v)\notin\wtilde V}}\ (\epsilon n)^{k\cdot\deg_{H}(v)/2}\prod_{u\in V(H)\setminus\{v\}}\deg\bigl(g(u)\bigr)^{k\cdot\deg_{H}(u)/2}\\
                                                     &\leq{n\over\abs{\Aut H}}\sum_{v\in V(H)}(\epsilon n)^{k\cdot\deg_{H}(v)/2}\sum_{g\colon V(H-v)\to V}\ \prod_{u\in V(H-v)}\deg\bigl(g(u)\bigr)^{k\cdot\deg_H(u)/2}.
        \end{align*}
        From here, we use the fact that $k\geq 2$ and proceed by the same steps in \Cref{prop:easyupper} to bound
        \begin{align*}
            (k!)^m\numb(G,\bup{H}{k})-\wtilde M &\leq{n\over\abs{\Aut H}}\sum_{v\in V(H)}(\epsilon n)^{k\cdot\deg_H(v)/2}\prod_{u\in V(H-v)}\biggl(\sum_{x\in V}\deg(x)^{k\cdot\deg_H(u)/2}\biggr)\\
                                                     &\leq{n\over\abs{\Aut H}}\sum_{v\in V(H)}(\epsilon n)^{k\cdot\deg_H(v)/2}\prod_{u\in V(H-v)}\biggl(\sum_{x\in V}\deg(x)\biggr)^{k\cdot\deg_H(u)/2}\\
                                                     &\leq{n\over\abs{\Aut H}}\sum_{v\in V(H)}(\epsilon n)^{k\cdot\deg_H(v)/2}\cdot(2Cn)^{km-k\cdot\deg_H(v)/2}\\
                                                     &\leq {n\over\abs{\Aut H}}\cdot\abs{V(H)}\cdot\epsilon^{k/2}\cdot (2C)^{km}\cdot n^{km}=O(\epsilon^{k/2} n^{km+1}).\qedhere
        \end{align*}
    \end{proof}

    We turn our attention now to bounding $\wtilde M$.
    Set $S\eqdef\sum_{uv\in{\wtilde V\choose 2}}\deg(u,v)$ and define the probability mass $\mu$ on ${\wtilde V\choose 2}$ by $\mu(uv)=\deg(u,v)/S$.
    By applying \Cref{prop:codegreebound}, we see that
    \[
        \wtilde M =\optb(\mu;H,k)\cdot S^{km}\leq\optb(H,k)\cdot\bigl(n+O(\epsilon^{-4})\bigr)^{km}=\optb(H,k)\cdot n^{km}\cdot\biggl(1+O\biggl({1\over n\epsilon^4}\biggr)\biggr)^{km}
    \]
    Therefore, if $\epsilon^4 n\to\infty$, we have
    \begin{equation}\label{eqn:epsntoinfity}
        \wtilde M\leq\optb(H,k)\cdot n^{km}+O\biggl({n^{km-1}\over\epsilon^4}\biggr).
    \end{equation}

    From here, we break into cases to conclude the proof.
    \medskip

    \textbf{Case: $k\cdot\bigl(\delta(H)-1\bigr)\geq 2$.}
    Select $\epsilon=n^{-1/(k+4)}$.
    Since $k\geq 1$, we have $\epsilon^4 n\to\infty$; hence we can apply \cref{eqn:epsntoinfity} and \Cref{claim:mindeg2} to bound
    \begin{align*}
        \numb(G,\bup{H}{k}) &\leq {\wtilde M\over (k!)^m}+O(\epsilon^k n^{km})\leq {\optb(H,k)\over (k!)^m}\cdot n^{km}+O\biggl({n^{km-1}\over\epsilon^4}\biggr)+O(\epsilon^k n^{km})\\
                      &={\optb(H,k)\over (k!)^m}\cdot n^{km}+O(n^{km- k/(k+4)}).
    \end{align*}
    \medskip

    \textbf{Case: $\delta(H)=1$ and $k\geq 9$.}
    Select $\epsilon=n^{-4/(k+8)}$.
    Since $k\geq 9$, we have $\epsilon^4 n\to\infty$; hence we can apply \cref{eqn:epsntoinfity} and \Cref{claim:mindeg1} to bound
    \begin{align*}
        \numb(G,\bup{H}{k}) &\leq {\wtilde M\over (k!)^m}+O(\epsilon^{k/2} n^{km+1})\leq {\optb(H,k)\over (k!)^m}\cdot n^{km}+O\biggl({n^{km-1}\over\epsilon^4}\biggr)+O(\epsilon^{k/2} n^{km+1})\\
                      &={\optb(H,k)\over (k!)^m}\cdot n^{km}+O(n^{km-1+16/(k+8)}).\qedhere
    \end{align*}
\end{proof}

Before moving on, we make a couple remarks.
\begin{remark}
    It is not difficult to argue that if $H$ is a graph on $m$ edges with no isolated vertices, $k$ is a positive integer and $C\geq 2$, then
    \[
        \numb_{\gcl C}(n,\bup{H}{k})\geq{\optb(H,k)\over(k!)^m}\cdot n^{km}-o(n^{km}),
    \]
    so \Cref{lem:tumorb} is asymptotically tight.
    Indeed, fix a finite set $X$ and a probability mass $\mu$ on ${X\choose 2}$.
    For a sufficiently large integer $n$, let $G$ be the edge-blow-up of $K_X$ formed by blowing up each edge $e\in{X\choose 2}$ into a set of size $\lfloor n\cdot\mu(e)\rfloor$.
    Then, one can show that $G\in\gcl 2$ and
    \[
        \numb(G,\bup{H}{k})\geq{\optb(\mu;H,k)\over(k!)^m}\cdot n^{km}-O(n^{km-1}).
    \]
\end{remark}
\begin{remark}\label{rem:planarblowup}
    For a finite set $X$ and a probability mass $\mu$ on ${X\choose 2}$, let $\sg\mu$ be the graph with vertex-set $X$ and edge-set $\supp\mu$.
    By following the proof of \Cref{lem:tumorp} more diligently, one can show that if $H$ is planar and $k\cdot\bigl(\delta(H)-1\bigr)\geq 2$, then
    \[
        \numb_\plan(n,\bup{H}{k})={\optb_\plan(H,k)\over(k!)^m}\cdot n^{km}+o(n^{km}),
    \]
    where
    \[
        \optb_\plan(H,k)=\sup\bigl\{\optb(\mu;H,k):\sg\mu\in\plan\bigr\}.
    \]
    Despite our beliefs when it comes to this same refinement in \Cref{lem:tumorp} (see~\Cref{rem:planarpath}), this could actually be an important refinement for certain planar graphs $H$.
    For instance, we believe that $\optb(K_5^-,1)>\optb_\plan(K_5^-,1)$ where $K_5^-$ is the $5$-clique minus an edge (we discuss this further in \Cref{sec:largeblowup}).

    In any case, we do not know how to prove a similar refinement in the case that $\delta(H)=1$ and $k\geq 9$.
\end{remark}

\paragraph{Reduction lemma for even cycles.}
In order to prove \Cref{lem:tumorc}, we will first need a straight-forward upper bound on the number of even paths in a graph.
\begin{prop}\label{prop:oddpath}
    If $G=(V,E)$ is any graph and $m$ is a positive integer, then
    \[
        \numb(G,P_{2m})\leq {(2\abs{E})^{m}\over 2}.
    \]
\end{prop}
\begin{proof}
    Label $V(P_{2m})=\{p_1,p_2,\dots,p_{2m}\}$ in consecutive order.
    For $(v_1,\dots,v_m)\in(V)_m$, consider the copies of $P_{2m}$ in $G$ wherein $v_i$ plays the role of $p_{2i}$.
    Observe that there are then at most $\deg(v_1)$ choices for the image of $p_1$ and at most $\deg(v_i,v_{i+1})$ choices for the image of $p_{2i+1}$ for all $i\in[m-1]$.
    Since $\deg(v_i,v_{i+1})\leq\deg(v_{i+1})$, we can therefore bound
    \begin{align*}
        \numb(G,P_{2m-1}) &\leq {1\over 2}\sum_{\mbf v\in(V)_m}\deg(v_1)\biggl(\prod_{i=1}^{m-1}\deg(v_i,v_{i+1})\biggr)\leq{1\over 2}\sum_{\mbf v\in(V)_m}\prod_{i=1}^m\deg(v_i)\\
                          &\leq {1\over 2}\sum_{v_1,\dots,v_{m}\in V}\prod_{i=1}^{m}\deg(v_i)={1\over 2}\biggl(\sum_{v\in V}\deg(v)\biggr)^{m}= {(2\abs{E})^{m}\over 2}.\qedhere
    \end{align*}
\end{proof}

We require additionally a simple observation about $2$-colorings of $C_m$.
\begin{prop}\label{prop:2color}
    Fix $m\geq 2$.
    For any $2$-coloring $\chi\colon\Z/m\Z\to\{0,1\}$, there is some $i\in\Z/m\Z$ for which either $\chi(i)=\chi(i+2)=0$ or $\chi(i)=\chi(i+3)=1$.
\end{prop}
\begin{proof}
    Suppose for the sake of contradiction that the claim does not hold.
    Since we are done if $\chi\equiv 1$, we may suppose, without loss of generality, that $\chi(0)=0$.
    This then implies that $\chi(-2)=\chi(2)=1$.
    But then $\chi(1)=\chi(-1)=0$; a contradiction.
\end{proof}

We are now ready to prove the reduction lemma for even cycles.

Recalling \cref{eqn:codegupper} from the proof of \Cref{lem:tumorb}, we know that for a graph $G$,
\begin{align*}
    \numb(G,C_4) &\leq\sum_{uv\in{V\choose 2}}{\deg(u,v)\choose 2},\quad\text{and}\\
    \numb(G,C_{2m}) &\leq\sum_{H\in\cp(V,C_m)}\ \prod_{xy\in E(H)}\deg(x,y)\quad\text{for }m\geq 3.
\end{align*}
We will not use either of these inequalities directly, but it will be helpful to keep them in mind throughout the following proof.

\begin{proof}[Proof of \Cref{lem:tumorc}]
    Let $G=(V,E)\in\gcl C$ be a graph on $n$ vertices.
    Fix $\epsilon=\epsilon(n)>0$ to be chosen later and define
    \[
        \wtilde V\eqdef\bigl\{v\in V:\deg(v)\geq\epsilon n\bigr\}.
    \]
    We denote an element of $\cp(G,C_{2m})$ by a tuple $(u_1,\dots,u_{2m})$, which is a list of the vertices of the cycle in some cyclic order.
    We define the following sets
    \begin{align*}
        \tp{Good} &\eqdef \bigl\{(v_1,\dots,v_{2m})\in\cp(G,C_{2m}): v_1,v_3,\dots,v_{2m-1}\in\wtilde V\text{ or }v_2,v_4,\dots,v_{2m}\in\wtilde V\bigr\},\\
        \tp{Bad} &\eqdef \cp(G,C_{2m})\setminus\tp{Good},\\
        \tp{Big} &\eqdef \bigl\{(v_1,\dots,v_{2m})\in\tp{Bad}: v_i,v_{i+2}\in\wtilde V\text{ for some }i\in[2m]\bigr\},\\
        \tp{Small} &\eqdef \bigl\{(v_1,\dots,v_{2m})\in\tp{Bad}: v_i,v_{i+3}\notin\wtilde V\text{ for some }i\in[2m]\bigr\}.
    \end{align*}
    Thanks to \Cref{prop:2color}, we know that $\tp{Bad}=\tp{Big}\cup\tp{Small}$.
    We aim to show that $\numb(G,C_{2m})\approx\abs{\tp{Good}}$.
    To do so, we must show that both $\tp{Big}$ and $\tp{Small}$ are both of insignificant size.
    \begin{claim}\label{claim:bigbound}
        $\abs{\tp{Big}}\leq O(\epsilon n^m)+O(n^{m-1}/\epsilon^3)$.
    \end{claim}
    \begin{proof}
        If $m=2$, then $\tp{Big}=\varnothing$ and so the claim holds.
        Hence, we may suppose that $m\geq 3$.
        Fix $H=(u_1,\dots,u_{2m})\in\tp{Big}$; without loss of generality, we may suppose that $u_1,u_3\in\wtilde V$.
        Since $H\in\tp{Bad}$, there must be some $i\in\{5,7,\dots,2m-1\}$ for which $\deg(u_i)<\epsilon n$, and so we bound
        \[
            \prod_{i=1}^m\deg(u_{2i-1},u_{2i+1})\leq\epsilon n\cdot \deg(u_1,u_3)\cdot\prod_{i=3}^m\deg(u_{2i-1}).
        \]
        By appealing additionally to \Cref{prop:codegreebound}, we can therefore crudely bound
        \begin{align*}
            \abs{\tp{Big}} &\leq \sum_{\substack{\mbf v\in (V)_m:\\ v_1,v_2\in\wtilde V,\\ v_i\notin\wtilde V\text{ for some }i\in[3,m]}}\prod_{i=1}^m\deg(v_i,v_{i+1})\leq \sum_{\substack{v_1\neq v_2\in\wtilde V,\\ v_3,\dots,v_m\in V}}\epsilon n\cdot\deg(v_1,v_2)\cdot\prod_{i=3}^m\deg(v_i)\\
                           &=\epsilon n\cdot\biggl(\sum_{v_1\neq v_2\in\wtilde V}\deg(v_1,v_2)\biggr)\biggl(\sum_{v\in V}\deg(v)\biggr)^{m-2}\leq\epsilon n\cdot\biggl(2n+O\biggl({1\over\epsilon^4}\biggl)\biggr)\cdot (2Cn)^{m-2}\\
                           &=O(\epsilon n^m)+O\biggl({n^{m-1}\over\epsilon^3}\biggr).\qedhere
        \end{align*}
    \end{proof}
    \begin{claim}\label{claim:smallbound}
        $\abs{\tp{Small}}\leq O(\epsilon n^m)$.
    \end{claim}
    \begin{proof}
        Fix $(u_1,\dots,u_{2m})\in\tp{Small}$; without loss of generality, we may suppose that $u_{2m-2},u_1\notin\wtilde V$.
        Observe that $u_1,\dots,u_{2m-2}$ forms a copy of $P_{2m-2}$ and that the edge $u_{2m-1}u_{2m}$ has both end-points in $N(u_1)\cup N(u_{2m-2})$.
        Therefore, by applying \Cref{prop:oddpath} and using the fact that $G\in\gcl C$, we see that
        \begin{align*}
            \abs{\tp{Small}} &\leq 2\cdot\numb(G,P_{2m-2})\cdot\max_{u\neq v\in V\setminus\wtilde V}\abss{E\bigl(G[N(u)\cup N(v)]\bigr)}\\
                             &\leq (2\abs{E})^{m-1}\cdot \max_{u\neq v\in V\setminus\wtilde V}C\cdot\abss{N(u)\cup N(v)}\\
                             &\leq (2Cn)^{m-1}\cdot C\cdot 2\epsilon n =O(\epsilon n^m).\qedhere
        \end{align*}
    \end{proof}
    We now deal with $\tp{Good}$.
    Define $S\eqdef\sum_{uv\in{\wtilde V\choose 2}}\deg(u,v)$ and let $\mu$ be the probability mass on ${\wtilde V\choose 2}$ defined by $\mu(uv)=\deg(u,v)/S$.
    If $m=2$, then we can bound
    \begin{align*}
        \abs{\tp{Good}} &\leq\sum_{uv\in {\wtilde V\choose 2}}{\deg(u,v)\choose 2}\leq {1\over 2}\sum_{uv\in{\wtilde V\choose 2}}\deg(u,v)^2={S^2\over 2}\sum_{uv\in{\wtilde V\choose 2}}\mu(uv)^2\\
                        &={S^2\over 2}\sum_{H\in\cp(\wtilde V,K_2)}\mu(H)^2={\optb(\mu;K_2,2)\over 2}\cdot S^2\leq{\optb(K_2,2)\over 2}\cdot S^2.
    \end{align*}
    Similarly, if $m\geq 3$, then we can bound
    \begin{align*}
        \abs{\tp{Good}} &\leq\sum_{H\in\cp(\wtilde V,C_m)}\ \prod_{xy\in E(H)}\deg(x,y)=S^m\sum_{H\in\cp(\wtilde V,C_m)}\ \prod_{xy\in E(H)}\mu(xy)\\
                        &=S^m\sum_{H\in\cp(\wtilde V,C_m)}\mu(H)=\optb(\mu;C_m,1)\cdot S^m\leq\optb(C_m,1)\cdot S^m.
    \end{align*}
    Thus, by applying \Cref{prop:codegreebound} and setting $B_2=\optb(K_2,2)/2$ and $B_m=\optb(C_m,1)$ for all $m\geq 3$, we have shown that
    \[
        \abs{\tp{Good}}\leq B_m\cdot S^m\leq B_m\bigl(n+O(1/\epsilon^4)\bigr)^m= B_m\cdot n^m+O\biggl({n^{m-1}\over\epsilon^4}\biggr),
    \]
    provided $\epsilon^4 n\to\infty$.
    Therefore, by selecting $\epsilon=n^{-1/5}$ and applying \Cref{claim:bigbound,claim:smallbound}, we bound
    \begin{align*}
        \numb(G,C_{2m}) &=\abs{\tp{Good}}+\abs{\tp{Bad}}\leq \abs{\tp{Good}}+\abs{\tp{Big}}+\abs{\tp{Small}}\\
                        &\leq B_m\cdot n^m+O\biggl({n^{m-1}\over\epsilon^4}\biggr)+O(\epsilon n^m)+O\biggl({n^{m-1}\over\epsilon^3}\biggr) =B_m\cdot n^m+O(n^{m-1/5}).\qedhere
    \end{align*}
\end{proof}

\section{Odd paths}\label{sec:paths}
Thanks to \Cref{lem:tumorp}, in order to bound $\numb_\plan(G,P_{2m+1})$ from above, it suffices to find upper bounds on $\optp(m)$.
Recall that for a finite set $X$ and a probability mass $\mu$ on ${X\choose 2}$,
\[
    \optp(\mu;m)=\sum_{\mbf x\in (X)_m}\bar\mu(x_1)\biggl(\prod_{i=1}^{m-1}\mu(x_ix_{i+1})\biggr)\bar\mu(x_m),
\]
where $\bar\mu(x)=\sum_{y\in X\setminus\{x\}}\mu(xy)$.
\medskip

First, we handle the case of $m=2$.
\begin{prop}\label{prop:p4}
    $\optp(2)=2$.
\end{prop}
\begin{proof}
    The lower bound is realized if $\abs{\supp\mu}=1$.
    \medskip

    For the upper bound, fix a finite set $X$ and a probability mass $\mu$ on ${X\choose 2}$.
    Define the matrix $M\in\R^{X\times X}$ by $M_{xy}=\mu(xy)$ under the convention that $\mu(xx)=0$.
    Observe that $M$ is a symmetric, non-negative matrix all of whose row-sums are bounded above by $1$.
    In particular, the largest eigenvalue of $M$ is at most $1$ (c.f.~\cite[Lemma 8.1.21]{horn_matrix}).
    Thus, by applying standard facts about the Rayleigh quotient (c.f.~\cite[Theorem 4.2.2]{horn_matrix}) and using the fact that $\sum_{x\in X}\bar\mu(x)=2$, we bound
    \[
        \optp(\mu;2)=\sum_{x\neq y\in X}\bar\mu(x)\mu(xy)\bar\mu(y)=\langle\bar\mu, M\bar\mu\rangle\leq\langle\bar\mu,\bar\mu\rangle=\sum_{x\in X}\bar\mu(x)^2\leq\sum_{x\in X}\bar\mu(x)=2.\qedhere
    \]
\end{proof}
From here, we have a quick proof of the asymptotic result of Ghosh et al.~\cite{ghosh_planarp5}, albeit with a worse error term.
\begin{proof}[Proof of \cref{known:p5} of \Cref{thm:alreadyknown}]
    By applying \Cref{lem:tumorp} and \Cref{prop:p4}, we bound
    \[
        \numb_\plan(n,P_5)\leq\numb_{\gcl 3}(n,P_5)\leq{\optp(2)\over 2}\cdot n^3+O(n^{14/5})=n^3+O(n^{14/5}).\qedhere
    \]
\end{proof}

Next, we establish a general upper bound on $\optp(m)$.

\begin{theorem}\label{thm:probp}
    For any $m\geq 3$,
    \[
        \optp(m)\leq{1\over(m-1)!}.
    \]
\end{theorem}
\begin{proof}
    Fix a finite set $X$ and a probability mass $\mu$ on ${X\choose 2}$.
    The key to this bound is to interpret $\optp(\mu;m)$ as the probability of some event in a probability space defined by $\mu$.
    Intuitively $\optp(\mu;m)$ is the probability that if we independently sample edges $e_1,\dots,e_{m+1}$ from $\mu$, then $e_2,\dots,e_m$ form a path with vertices $x_1,\dots,x_m$, $e_1$ is incident to $x_1$ and $e_{m+1}$ is incident to $x_m$.
    We now make this intuition precise.
    \medskip

    For a tuple $\mbf x\in (X)_m$, define the sets
    \begin{align*}
        \mcal E(\mbf x) &\eqdef\biggl\{(e_2,\ldots,e_m)\in{X\choose 2}^{m-1}:\{e_2,\ldots,e_m\}=\{x_1x_2,x_2x_3,\ldots,x_{m-1}x_m\}\biggr\},\qquad\text{and}\\
        \mcal L(\mbf x) &\eqdef\biggl\{(e_1,e_{m+1})\in{X\choose 2}^2:e_1\ni x_1\text{ and }e_{m+1}\ni x_m\biggr\}.
    \end{align*}
    Observe that $\Pr_{\mu^2}[\mcal L(\mbf x)]=\bar\mu(x_1)\cdot\bar\mu(x_m)$ and that
    \[
        \mu^{m-1}(\mbf e)=\prod_{i=1}^{m-1}\mu(x_ix_{i+1}),\qquad\text{for all }\mbf e\in\mcal E(\mbf x),
    \]
    where $\mu^j$ is the product distribution induced on ${X\choose 2}^j$ by $\mu$.

    We can therefore write
    \begin{align*}
        \optp(\mu;m) &= \sum_{\mbf x\in (X)_m}\bar\mu(x_1)\biggl(\prod_{i=1}^{m-1}\mu(x_ix_{i+1})\biggr)\bar\mu(x_m)=\sum_{\mbf x\in(X)_m}\biggl({1\over\abs{\mcal E(\mbf x)}}\sum_{\mbf e\in\mcal E(\mbf x)}\mu^{m-1}(\mbf e)\biggr)\cdot\Pr_{\mu^2}[\mcal L(\mbf x)]\\
                     &= {1\over(m-1)!}\sum_{\mbf x\in(X)_m}\Pr_{\mu^{m-1}}[\mcal E(\mbf x)]\cdot\Pr_{\mu^2}[\mcal L(\mbf x)]={1\over(m-1)!}\sum_{\mbf x\in(X)_m}\Pr_{\mu^{m+1}}[\mcal E(\mbf x)\times\mcal L(\mbf x)].
    \end{align*}

    For $\mbf x\in(X)_m$, consider the reverse tuple $\wtilde{\mbf x}\in(X)_m$ where $\wtilde x_i=x_{m+1-i}$. The events $\mcal E(\mbf x)\times\mcal L(\mbf x)$ and $\mcal E(\mbf y)\times\mcal L(\mbf y)$ are almost always disjoint when $\mbf x\neq\mbf y$. The only circumstance in which they double-count the same event is when $\mbf y=\wtilde{\mbf x}$, in which case $\mcal E(\mbf x)=\mcal E(\wtilde{\mbf x})$ and $\mcal L(\mbf x)\cap\mcal L(\wtilde{\mbf x})=\{(x_1x_m, x_1x_m)\}$.
    Indeed, for $\mbf x,\mbf y\in(X)_m$, we have $\{x_1x_2,\dots,x_{m-1}x_m\}=\{y_1y_2,\dots,y_{m-1}y_m\}$ if and only if $\mbf y\in\{\mbf x,\wtilde{\mbf x}\}$; thus
    \[
        \mcal E(\mbf x)\cap\mcal E(\mbf y)=\begin{cases}
            \mcal E(\mbf x) & \text{if }\mbf y\in\{\mbf x,\wtilde{\mbf x}\},\\
            \varnothing & \text{otherwise},
        \end{cases}
    \]
    and
    \[
        \bigl(\mcal E(\mbf x)\times\mcal L(\mbf x)\bigr)\cap\bigl(\mcal E(\mbf y)\times\mcal L(\mbf y)\bigr)=\begin{cases}
            \mcal E(\mbf x)\times\mcal L(\mbf x) & \text{if }\mbf y=\mbf x,\\
            \mcal E(\mbf x)\times\{(x_1x_m,\ x_1x_m)\} & \text{if }\mbf y=\wtilde{\mbf x},\\
            \varnothing & \text{otherwise}.
        \end{cases}
    \]
    Therefore, by grouping together $\mbf x$ and $\wtilde{\mbf x}$, we compute
    \begin{align*}
        (m-1)!\cdot\optp(\mu;m) &= \sum_{\mbf x\in (X)_m}\Pr_{\mu^{m+1}}[\mcal E(\mbf x)\times\mcal L(\mbf x)]\\
                                &=\Pr_{\mu^{m+1}}\biggl[\bigcup_{\mbf x\in (X)_m}\bigl(\mcal E(\mbf x)\times\mcal L(\mbf x)\bigr)\biggr]+\Pr_{\mu^{m+1}}\biggl[\bigcup_{\mbf x\in (X)_m}\bigl(\mcal E(\mbf x)\times\{(x_1x_m,x_1x_m)\}\bigr)\biggr].
    \end{align*}

    Next, by writing
    \[
        \mcal E\eqdef\bigcup_{\mbf x\in(X)_m}\mcal E(\mbf x),
    \]
    we bound
    \begin{align*}
        \Pr_{\mu^{m+1}}\biggl[\bigcup_{\mbf x\in (X)_m}\bigl(\mcal E(\mbf x)\times\mcal L(\mbf x)\bigr)\biggr] &\leq \Pr_{\mu^{m+1}}\biggl[\mcal E\times\bigcup_{\mbf x\in(X)_m}\mcal L(\mbf x)\biggr]\leq\Pr_{\mu^{m-1}}[\mcal E],
    \end{align*}
    and
    \begin{align*}
        \Pr_{\mu^{m+1}}\biggl[\bigcup_{\mbf x\in (X)_m}\bigl(\mcal E(\mbf x)\times\{(x_1x_m,x_1x_m)\}\bigr)\biggr] &\leq\Pr_{\mu^{m+1}}\biggl[\mcal E\times\bigcup_{\mbf x\in(X)_m}\{(x_1x_m,x_1x_m)\}\biggr]\\
                                                                                                                          &=\Pr_{\mu^{m-1}}[\mcal E]\cdot\Pr_{\mu^2}\biggl[\bigcup_{e\in{X\choose 2}}\{(e,e)\}\biggr] =\Pr_{\mu^{m-1}}[\mcal E]\cdot\sum_{e\in{X\choose 2}}\mu(e)^2.
    \end{align*}
    Any member of $\mcal E$ has the property that its coordinates are distinct members of ${X\choose 2}$; hence, since $m\geq 3$, we can bound
    \begin{align*}
        \Pr_{\mu^{m-1}}[\mcal E] &\leq\Pr_{\mu^{m-1}}\biggl[\biggl\{(e_2,e_3,\ldots,e_m)\in{X\choose 2}^{m-1}:e_2,e_3,\dots,e_m\text{ distinct}\biggr\}\biggr]\\
                                 &\leq\Pr_{\mu^2}\biggl[\biggl\{(e_2,e_3)\in{X\choose 2}^2:e_2\neq e_3\biggr\}\biggr]=1-\sum_{e\in{X\choose 2}}\mu(e)^2.
    \end{align*}
    Putting everything together, we have shown that
    \begin{align*}
        (m-1)!\cdot\optp(\mu;m) &\leq\Pr_{\mu^{m-1}}[\mcal E]+\Pr_{\mu^{m-1}}[\mcal E]\cdot\sum_{e\in{X\choose 2}}\mu(e)^2\leq\biggl(1-\sum_{e\in{X\choose 2}}\mu(e)^2\biggr)\biggl(1+\sum_{e\in{X\choose 2}}\mu(e)^2\biggr)\\
                                &=1-\biggl(\sum_{e\in{X\choose 2}}\mu(e)^2\biggr)^2\leq 1,
    \end{align*}
    which establishes the claim.
\end{proof}

\medskip


\subsection{Paths of order 7}
In this section, we prove \Cref{thm:oddpath}.
The main content in this section is the proof that $\optp(3)=8/27$, which hinges on the following general inequality.
We note that the following lemma is a special case of the much more general \Cref{thm:largek}, but we give a direct and self-contained proof here.
\begin{lemma}\label{lem:aequalb}
    If $a_1,\dots,a_n\geq 0$, then
    \[
        \biggl(\sum_i a_i^2\biggr)^2-\sum_i a_i^4\leq {1\over 8}\biggl(\sum_i a_i\biggr)^4.
    \]
\end{lemma}
\begin{proof}
    We notice first that the claim is trivial if $a_i=0$ for all $i$.
    Furthermore, scaling the $a_i$'s by any positive constant leaves the inequality invariant.
    As such, we may suppose that $\sum_i a_i=1$.

    Therefore, noting that $\bigl(\sum_i x_i^2\bigr)^2-\sum_i x_i^4=\sum_{i\neq j}x_i^2x_j^2$, it suffices to show that
    \begin{equation}\label{eqn:optaequalb}
        \begin{array}{cl}
            \max & \sum_{i\neq j}x_i^2x_j^2\\
            \text{s.t.} & \sum_i x_i=1\\
                        & x_i\geq 0 \qquad\text{for all }i\in[n],
        \end{array}
    \end{equation}
    is bounded above by $1/8$.
    Let $a_1,\dots,a_n$ denote an optimal solution to \cref{eqn:optaequalb}; without loss of generality, we may suppose that $a_1\geq\dots\geq a_n>0$.
    Additionally, let $M$ denote the optimal value, that is, $M=\sum_{i\neq j}a_i^2a_j^2$.
    We may certainly suppose that $n\geq 2$ since otherwise $M=0$.

    By applying the KKT conditions (\Cref{KKT}) to \cref{eqn:optaequalb}, we find that there is some fixed $\lambda\in\R$ for which
    \begin{equation}\label{eqn:regularity}
        a_i\sum_{j:\ j\neq i}a_j^2 = \lambda\qquad\text{for all }i\in[n].
    \end{equation}
    From here, we use the fact that $\sum_ia_i=1$ to determine,
    \begin{equation}\label{eqn:lambdavalue}
        \lambda = \sum_i a_i\lambda = \sum_ia_i^2\sum_{j:\ j\neq i}a_j^2=\sum_{i\neq j}a_i^2a_j^2=M.
    \end{equation}
    \medskip

    Now, consider the numbers $b_1,\dots,b_{n-1}$ defined by
    \[
        b_i=a_i/(1-a_n),
    \]
    which are well-defined since $n\geq 2$ and hence $a_n<1$.
    Note that $b_i>0$ and $\sum_ib_i=1$.
    Therefore,
    \begin{align*}
        M \geq \sum_{i\neq j}b_i^2b_j^2 &={1\over(1-a_n)^4}\sum_{\substack{i,j\in[n-1]:\\ i\neq j}}a_i^2a_j^2={1\over(1-a_n)^4}\biggl(\sum_{i\neq j}a_i^2a_j^2-2a_n^2\sum_{j=1}^{n-1}a_j^2\biggr)\\
                                        &={1\over (1-a_n)^4}\bigl(M-2Ma_n\bigr)=M\cdot{1-2a_n\over(1-a_n)^4},
    \end{align*}
    where the penultimate equality follows from \cref{eqn:regularity,eqn:lambdavalue}.
    We conclude that $1-2a_n\leq(1-a_n)^4$ and thus $a_n\geq 0.45$.
    Since $a_1\geq\dots\geq a_n$, this then implies that $n=2$.
    Thus, we apply the AM--GM inequality to finally bound
    \[
        M=2a_1^2a_2^2\leq2\biggl({a_1+a_2\over 2}\biggr)^4={1\over 8}.\qedhere
    \]
\end{proof}

We will apply the following direct corollary of \Cref{lem:aequalb}.

\begin{corollary}\label{cor:offdiag}
    If $a_1,\dots,a_n,b_1,\dots,b_n\geq 0$, then
    \[
        \biggl(\sum_i a_ib_i\biggr)^2-\sum_i a_i^2b_i^2\leq{1\over 8}\biggl(\sum_i a_i\biggr)^2\biggl(\sum_i b_i\biggr)^2.
    \]
\end{corollary}
\begin{proof}
    By applying \Cref{lem:aequalb} and the Cauchy--Schwarz inequality, we bound
    \begin{align*}
        \biggl(\sum_i a_ib_i\biggr)^2-\sum_i a_i^2b_i^2 &= \biggl(\sum_i (\sqrt{a_ib_i})^2\biggr)^2-\sum_i (\sqrt{a_ib_i})^4\\
                                                        &\leq {1\over 8}\biggl(\sum_i \sqrt{a_ib_i}\biggr)^4\leq{1\over 8}\biggl(\sum_i a_i\biggr)^2\biggl(\sum_i b_i\biggr)^2.\qedhere
    \end{align*}
\end{proof}

We can now determine $\optp(3)$.
\begin{theorem}\label{thm:p6}
    $\optp(3)=8/27$.
\end{theorem}
\begin{proof}
    To prove the lower bound, let $\mu$ be the uniform distribution on ${[3]\choose 2}$.
    Then
    \[
        \optp(\mu;3)=\sum_{(x,y,z)\in([3])_3}\bar\mu(x)\mu(xy)\mu(yz)\bar\mu(z)=3!\cdot{2\over 3}\cdot{1\over 3}\cdot{1\over 3}\cdot{2\over 3}={8\over 27}.
    \]

    To establish the upper bound, fix a finite set $X$ and let $\mu$ be any probability distribution on ${X\choose 2}$.
    We begin by writing
    \begin{align*}
        \optp(\mu;3) &=\sum_{(x,y,z)\in(X)_3}\bar\mu(x)\mu(xy)\mu(yz)\bar\mu(z)=\sum_{y\in X}\sum_{x\in X\setminus\{y\}}\bar\mu(x)\mu(xy)\sum_{z\in X\setminus\{x,y\}}\bar\mu(z)\mu(zy)\\
                     &=\sum_{y\in X}\biggl[\biggl(\sum_{x\in X\setminus\{y\}}\bar\mu(x)\mu(xy)\biggr)^2-\sum_{x\in X\setminus\{y\}}\bar\mu(x)^2\mu(xy)^2\biggr].
    \end{align*}
    Then, by applying \Cref{cor:offdiag} to the inner expression and using the fact that $\sum_{x\in X}\bar\mu(x)=2$, we bound
    \begin{align*}
        \optp(\mu;3) &\leq \sum_{y\in X}{1\over 8}\biggl(\sum_{x\in X\setminus\{y\}}\bar\mu(x)\biggr)^2\biggl(\sum_{x\in X\setminus\{y\}}\mu(xy)\biggr)^2 ={1\over 8}\sum_{y\in X}\bigl(2-\bar\mu(y)\bigr)^2\cdot\bar\mu(y)^2.
    \end{align*}
    We finally observe that the expression $x(2-x)^2$ for $0\leq x\leq 1$ is maximized when $x=2/3$, yielding a value of $32/27$.
    Therefore,
    \begin{align*}
        \optp(\mu;3) &\leq {1\over 8}\sum_{y\in X}\bar\mu(y)\cdot\bar\mu(y)\bigl(2-\bar\mu(y)\bigr)^2\leq {4\over 27}\sum_{y\in X}\bar\mu(y)={8\over 27}.\qedhere
    \end{align*}
\end{proof}

Now that we know $\optp(3)$, the proof of \Cref{thm:oddpath} follows quickly.
\begin{proof}[Proof of \Cref{thm:oddpath}]
    First, the graph $\bup{K_3}{\ell}$ where $\ell=\lfloor{n-3\over 3}\rfloor$ shows that
    \[
        \numb_\plan(n,P_7)\geq{4\over 27}\cdot n^4-O(n^3).
    \]

    Next, we apply \Cref{lem:tumorp} to bound
    \[
        \numb_\plan(n,P_{2m+1})\leq\numb_{\gcl 3}(n,P_{2m+1})\leq{\optp(m)\over 2}\cdot n^{m+1}+O(n^{m+4/5}),
    \]
    for all $m\geq 2$.
    Finally, \Cref{thm:p6} tells us that $\optp(3)=8/27$, and \Cref{thm:probp} tells us that $\optp(m)\leq 1/(m-1)!$ for all $m\geq 4$; hence the claim follows.
\end{proof}

\section{Edge-blow-ups and even cycles}\label{sec:blowup}
Thanks to \Cref{lem:tumorb,lem:tumorc}, in order to bound $\numb_\plan(n,\bup{H}{k})$ from above for various $H,k$, it suffices to prove upper bounds on $\optb(H,k)$.
Recall that for a finite set $X$ and a probability mass $\mu$ on ${X\choose 2}$,
\[
    \optb(\mu;H,k)=\sum_{H'\in\cp(X,H)}\mu(H')^k,
\]
where $\cp(X,H)$ is the set of copies of $H$ in $K_X$, and
\[
    \mu(H')=\prod_{e\in E(H')}\mu(e).
\]

We deal first with the case of $H=K_2$.
\begin{prop}\label{prop:k2}
    $\optb(K_2,k)=1$ for all $k\geq 1$.
\end{prop}
\begin{proof}
    The lower bound is realized if $\abs{\supp\mu}=1$.
    \medskip

    Let $X$ be a finite set and let $\mu$ be a probability mass on ${X\choose 2}$.
    Then
    \[
        \optb(\mu;K_2,k)=\sum_{K\in\cp(X,K_2)}\mu(K)^k=\sum_{e\in{X\choose 2}}\mu(e)^k\leq \biggl(\sum_{e\in {X\choose 2}}\mu(e)\biggr)^k=1.\qedhere
    \]
\end{proof}

Since $\plan\subseteq\gcl 3$, \cref{known:c4,known:k2k} of \Cref{thm:alreadyknown} follow immediately thanks to \Crefnosort{lem:tumorc,lem:tumorb}, respectively.
\medskip

Next, we prove a general upper bound on $\optb(H,k)$.
\begin{theorem}\label{thm:probb}
    If $H$ is a graph on $m$ edges with no isolated vertices and $k$ is a positive integer, then
    \[
        \optb(H,k)\leq{\ (k!)^m\over (km)!}.
    \]
\end{theorem}
\begin{proof}
    Fix a finite set $X$ and let $\mu$ be a probability mass on ${X\choose 2}$.
    The key to this bound is to relate $\optb(\mu;H,k)$ to an event in a probability space defined by $\mu$.
    Intuitively, $\optb(\mu;H,k)$ is the probability that $km$ edges sampled independently from $\mu$ form a copy of $H$ wherein each edge has multiplicity $k$.
    We now make this intuition precise.
    \medskip

    For $H'\in\cp(X,H)$, define the set
    \[
        \mcal C(H')\eqdef\biggl\{\mbf e\in{X\choose 2}^{km}: \text{each $e\in E(H')$ occurs exactly $k$ times in }\mbf e\biggr\}.
    \]
    Observe that
    \[
        \abs{\mcal C(H')}={km\choose k,\dots, k}={(km)!\over(k!)^m},
    \]
    and that $\mu(H')^k=\mu^{km}(\mbf e)$ for every $\mbf{e}\in\mcal C(H')$, where $\mu^{km}$ is the product distribution on ${X\choose 2}^{km}$ induced by $\mu$.

    Now, the events $\bigl\{\mcal C(H'):H'\in\cp(X,H)\bigr\}$ are pairwise disjoint since the entries of any $\mbf e\in\mcal C(H')$ uniquely define the edge-set of $H'$. Consequently, we can bound
    \begin{align*}
        \optb(\mu;H,k) &= \sum_{H'\in\cp(X,H)}\mu(H')^k = \sum_{H'\in\cp(X,H)}{1\over\abs{\mcal C(H')}}\sum_{\mbf e\in\mcal C(H')}\mu^{km}(\mbf e)\\
                       &={(k!)^m\over(km)!}\cdot\sum_{H'\in\cp(X,H)}\Pr_{\mu^{km}}[\mcal C(H')]={(k!)^m\over(km)!}\cdot\Pr_{\mu^{km}}\biggl[\bigcup_{H'\in\cp(X,H)}\mcal C(H')\biggr] \leq {(k!)^m\over(km)!}.\qedhere
    \end{align*}
\end{proof}

From here, we can immediately prove \Cref{thm:blowupbound}.
\begin{proof}[Proof of \Cref{thm:blowupbound}]
    First, \Cref{thm:probb} tells us that $\optb(H,k)\leq(k!)^m/(km)!$.
    Then, thanks to \Cref{lem:tumorb}, if $k\cdot\bigl(\delta(H)-1\bigr)\geq 2$ or if $\delta(H)=1$ and $k\geq 9$, then
    \[
        \numb_\plan(n,\bup{H}{k})\leq\numb_{\gcl 3}(n,\bup{H}{k})\leq{\optb(H,k)\over(k!)^m}\cdot n^{km}+o(n^{km})\leq{n^{km}\over (km)!}+o(n^{km}).\qedhere
    \]
\end{proof}

\subsection{The structure of optimal masses}
In this section, we establish structural properties about those masses which achieve $\optb(H,k)$, which will be used in the next sections in order to prove \Cref{thm:evencycle,thm:planarclique,thm:blowupexact}.

Of course, a priori, it is not even clear that $\optb(H,k)$ is ever achieved.
In fact, one can show that $\optb(K_{1,m},1)=1/m!$ for all $m\geq 2$, yet this value is never achieved.
Indeed, one can argue that for all $n\geq m\geq 2$,
\[
    \max\bigl\{\optb(\mu;K_{1,m},1):\abs{\supp\mu}\leq n\bigr\}={n\choose m}\cdot{1\over n^m}<{1\over m!}.
\]
The same phenomenon occurs for $\optb(mK_2,1)$ for $m\geq 2$ where $mK_2$ is the matching on $m$ edges. We conjecture that these are the only situations in which $\optb(H,k)$ is not achieved. See \Cref{cor:achieved} for partial results in this direction.

Despite this, for any fixed, finite set $X$ with at least two elements, the quantity $\max\bigl\{\optb(\mu;H,k):\supp\mu\subseteq{X\choose 2}\bigr\}$ exists, thanks to compactness.
\begin{defn}
    Let $H$ be a graph with no isolated vertices and let $k$ be a positive integer.
    For a finite set $X$, we denote by $\Opt(X;H,k)$ the set of all probability masses $\mu$ on ${X\choose 2}$ satisfying
    \[
    \optb(\mu;H,k)=\max\biggl\{\optb(\mu';H,k):\supp\mu'\subseteq{X\choose 2}\biggr\}.
    \]
    In the case that $\optb(H,k)$ is achieved, we denote by $\Opt(H,k)$ the set of all masses $\mu$ satisfying $\optb(\mu;H,k)=\optb(H,k)$.
\end{defn}

Fix a finite set $X$ and a probability mass $\mu$ on ${X\choose 2}$ and let $\sg\mu$ be the graph with vertex-set $X$ and edge-set $\supp\mu$.
Observe that $\optb(\mu;H,k)>0$ if and only if $\sg\mu$ has a copy of $H$; consequently, if $\optb(\mu;H,k)>0$, then $\abs{\supp\mu}\geq\abs{E(H)}$ and $\abs{\supp\bar\mu}\geq\abs{V(H)}$.
We see also that if $\abs X\geq\abs{V(H)}$ and $\mu\in\Opt(X;H,k)$, then $\sg\mu$ must contain a copy of $H$.
Additionally, we can determine such an optimal $\mu$ exactly if $\abs{\supp\mu}=\abs{E(H)}$:

\begin{prop}\label{prop:uniform}
    Let $H=(V,E)$ be a graph on $m$ edges with no isolated vertices and let $k$ be a positive integer.
    Fix any finite set $X$ with $\abs X\geq\abs V$ and fix $\mu\in\Opt(X;H,k)$.
    If $\abs{\supp\mu}=m$, then $\mu$ is the uniform distribution on $E(H')$ for some $H'\in\cp(X,H)$ and thus $\optb(\mu;H,k)=m^{-km}$.
\end{prop}
\begin{proof}
    We know that $\sg\mu$ contains a copy of $H$ since $\abs X\geq\abs V$ and $\mu\in\Opt(X;H,k)$.
    Since $\abs{\supp\mu}=m$, we conclude that $\sg\mu$ must in fact be a copy of $H$, possibly with isolated vertices.
    We can therefore apply the arithmetic--geometric mean inequality to bound
    \[
        \optb(\mu;H,k)=\prod_{e\in\supp\mu}\mu(e)^k\leq\biggl({1\over m}\sum_{e\in\supp\mu}\mu(e)\biggr)^{km}={1\over m^{km}},
    \]
    with equality if and only if $\mu(e)=1/m$ for every $e\in\supp\mu$.
\end{proof}

We next derive regularity conditions for the members of $\Opt(X;H,k)$.
\begin{lemma}\label{lem:regularity}
    Let $H$ be a graph on $m$ edges with no isolated vertices, let $k$ be a positive integer and fix a finite set $X$.
    If $\mu\in\Opt(X;H,k)$, then
    \begin{alignat*}{2}
        \mu(e)\cdot m\cdot\optb(\mu;H,k) &= \sum_{\substack{H'\in\cp(X,H):\\ E(H')\ni e}}\mu(H')^k && \text{for every }e\in{X\choose 2},\quad\text{and}\\
        \bar\mu(x)\cdot m\cdot\optb(\mu;H,k) &= \sum_{\substack{H'\in\cp(X,H):\\ V(H')\ni x}}\deg_{H'}(x)\cdot\mu(H')^k \qquad && \text{for every }x\in X.
    \end{alignat*}
\end{lemma}
\begin{proof}
    By the definition of $\beta$, we can write
    \[
        \begin{array}{ccl}
            \optb(\mu;H,k) = & \max & \sum_{H'\in\cp(X,H)}\ \prod_{e\in E(H')}x_e^k\\
                             &\text{s.t.} & \sum_{e\in {X\choose 2}}x_e=1\\
                             & & x_e\geq 0 \qquad\text{for all }e\in{X\choose 2}.
        \end{array}
    \]
    In particular, we can apply the KKT conditions (\Cref{KKT}) to $\mu$.
    By doing so, we find that there is some fixed $\lambda\in\R$ such that  $D(e)=\lambda$ for all $e\in\supp\mu$, where
    \[
        D(e)\eqdef\sum_{\substack{H'\in\cp(X,H):\\ E(H')\ni e}}\ \mu(e)^{k-1}\prod_{s\in E(H')\setminus\{e\}}\mu(s)^k.
    \]
    Of course, whether or not $e\in\supp\mu$, we always have
    \[
        \lambda\cdot\mu(e)=D(e)\cdot\mu(e)=\sum_{\substack{H'\in\cp(X,H):\\ E(H')\ni e}}\mu(H')^k.
    \]
    Using $\mbf 1[S]$ to denote the indicator function of an event $S$, we compute
    \begin{align*}
        \lambda &= \sum_{e\in{X\choose 2}}\lambda\cdot\mu(e)=\sum_{e\in{X\choose 2}}D(e)\cdot\mu(e) = \sum_{e\in{X\choose 2}}\ \sum_{\substack{H'\in\cp(X,H):\\ E(H')\ni e}}\mu(H')^k\\
                &= \sum_{H'\in\cp(X,H)}\mu(H')^k\cdot\sum_{e\in{X\choose 2}}\mbf 1[e\in E(H')]= m\cdot\optb(\mu;H,k),
    \end{align*}
    and so
    \[
        \mu(e)\cdot m\cdot\optb(\mu;H,k)=\mu(e)\cdot\lambda=\sum_{\substack{H'\in\cp(X,H):\\ E(H')\ni e}}\mu(H')^k,
    \]
    for every $e\in{X\choose 2}$.

    From here, we see also that for each $x\in X$,
    \begin{align*}
        \bar\mu(x)\cdot m\cdot\optb(\mu;H,k) &= \sum_{y\in X\setminus\{x\}}\mu(xy)\cdot m\cdot\optb(\mu;H,k)=\sum_{y\in X\setminus\{x\}}\sum_{\substack{H'\in\cp(X,H):\\ E(H')\ni xy}}\mu(H')^k\\
                                             &=\sum_{H'\in\cp(X,H)}\mu(H')^k\cdot\sum_{y\in X\setminus\{x\}}\mbf 1[xy\in E(H')]\\
                                             &=\sum_{\substack{H'\in\cp(X,H):\\ V(H')\ni x}}\deg_{H'}(x)\cdot\mu(H')^k.\qedhere
    \end{align*}
\end{proof}

These regularity conditions allow us to place bounds on the edge- and vertex-masses in an optimal mass.
\begin{lemma}\label{lem:massbound}
    Let $H$ be a graph on $m$ edges with no isolated vertices, let $k$ be a positive integer and fix a finite set $X$ with $\abs X\geq\abs{V(H)}$.
    If $\mu\in\Opt(X;H,k)$, then
    \begin{alignat*}{2}
        1-m\cdot\mu(e)&\leq\bigl(1-\mu(e)\bigr)^{km}\qquad && \text{for all }e\in\supp\mu,\quad\text{and}\\
        1-{m\over\delta(H)}\bar\mu(x) &\leq \bigl(1-\bar\mu(x)\bigr)^{km} && \text{for all }x\in\supp\bar\mu.
    \end{alignat*}
\end{lemma}
\begin{proof}
    Since $\abs X\geq\abs{V(H)}$, we know that $\optb(\mu;H,k)>0$.

    We prove first that $1-m\cdot\mu(e)\leq\bigl(1-\mu(e)\bigr)^{km}$ for any $e\in\supp\mu$.
    Fix any $e\in\supp\mu$.
    If $\mu(e)\geq 1/m$, then the claim is trivial; otherwise, $\mu(e)<1/m$, and we can define the mass $\mu'$ on ${X\choose 2}$ by
    \[
        \mu'(s)={1\over 1-\mu(e)}\cdot\begin{cases}
            0 & \text{if }s=e,\\
            \mu(s) & \text{otherwise}.
        \end{cases}
    \]
    Since $\mu\in\Opt(X;H,k)$, we apply \Cref{lem:regularity} to see that
    \begin{align*}
        \optb(\mu;H,k) \geq \optb(\mu';H,k) &={1\over(1-\mu(e))^{km}}\cdot\biggl(\optb(\mu;H,k)-\sum_{\substack{H'\in\cp(X,H):\\ E(H')\ni e}}\mu(H')^k\biggr)\\
                       &=\optb(\mu;H,k)\cdot{1-m\cdot\mu(e)\over (1-\mu(e))^{km}},
    \end{align*}
    which implies that $1-m\cdot\mu(e)\leq\bigl(1-\mu(e)\bigr)^{km}$.
    \medskip

    We prove next that $1-{m\over\delta}\bar\mu(x)\leq\bigl(1-\bar\mu(x)\bigr)^{km}$ for any $x\in\supp\bar\mu$, where $\delta=\delta(H)$.
    Fix any $x\in\supp\bar\mu$.
    If $\bar\mu(x)\geq \delta/m$, then the claim is trivial; otherwise, $\bar\mu(x)<\delta/m$, and we can define the mass $\mu'$ on ${X\choose 2}$ by
    \[
        \mu'(s)={1\over 1-\bar\mu(x)}\cdot\begin{cases}
            0 & \text{if }s\ni x,\\
            \mu(s) & \text{otherwise}.
        \end{cases}
    \]
    Since $\mu\in\Opt(X;H,k)$, we again apply \Cref{lem:regularity} to see that
    \begin{align*}
        \optb(\mu;H,k) \geq\optb(\mu';H,k) &={1\over(1-\bar\mu(x))^{km}}\cdot\biggl(\optb(\mu;H,k)-\sum_{\substack{H'\in\cp(X,H):\\ V(H')\ni x}}\mu(H')^k\biggr)\\
                       &\geq{1\over(1-\bar\mu(x))^{km}}\cdot\biggl(\optb(\mu;H,k)-{\bar\mu(x)\over\delta}\cdot m\cdot\optb(\mu;H,k)\biggr)\\
                       &=\optb(\mu;H,k)\cdot{1-{m\over\delta}\bar\mu(x)\over(1-\bar\mu(x))^{km}},
    \end{align*}
    which implies that $1-{m\over\delta}\bar\mu(x)\leq\bigl(1-\bar\mu(x)\bigr)^{km}$.
\end{proof}
We remark that one can show also that $\mu(e)\leq 1/m$ for all $e\in\supp\mu$ and that $\bar\mu(x)\leq\Delta(H)/m$ for all $x\in\supp\bar\mu$; however, we have not found any use for these inequalities.

\medskip

\Cref{lem:massbound} allows us to place lower-bounds on $\mu(e)$ for $e\in\supp\mu$ and on $\bar\mu(x)$ for $x\in\supp\bar\mu$ when $\mu$ is an optimal mass.
For instance, consider the inequality $1-m\cdot\mu(e)\leq\bigl(1-\mu(e)\bigr)^{km}$.
This inequality always holds if $k=1$, but if $k\geq 2$, then we observe that the curves $1-mx$ and $(1-x)^{km}$ intersect at $0$ and at a unique $x^*\in(0,1]$.
Furthermore, $1-mx>(1-x)^{km}$ for all $x\in(0,x^*)$ and $1-mx\leq(1-x)^{km}$ for all $x\in[x^*,1]$.
Therefore, if we can locate some $x\in(0,1]$ for which $1-mx>(1-x)^{km}$, then we will have shown that $\mu(e)>x$ for all $e\in\supp\mu$.
Similar reasoning can be applied to the inequality $1-{m\over\delta}\bar\mu(x)\leq\bigl(1-\bar\mu(x)\bigr)^{km}$; that is, if we can locate some $z\in(0,1]$ for which $1-{m\over\delta}z>(1-z)^{km}$, then we will have shown that $\bar\mu(x)>z$ for all $x\in\supp\bar\mu$.

Indeed, we will apply precisely this reasoning in order to establish \Cref{thm:optclique,thm:optc4,thm:largek}.
However, before we get to this, we first remark on a useful consequence of \Cref{lem:massbound}.
\begin{corollary}\label{cor:achieved}
    For a graph $H$ and a positive integer $k$, if $k\cdot\delta(H)\geq 2$, then $\optb(H,k)$ is achieved.
\end{corollary}
\begin{proof}
    Let $X$ be a finite set with $\abs X\geq\abs{V(H)}$ and fix any $\mu\in\Opt(X;H,k)$.
    By passing to a subset of $X$ if necessary, we may suppose that $\supp\bar\mu=X$.
    Thanks to compactness, in order to show that $\optb(H,k)$ is achieved, it suffices to show that $\abs X$ is bounded above by some constant depending only on $H$ and $k$.

    Set $\delta=\delta(H)$, $m=\abs{E(H)}$ and fix $x\in X$ with $\bar\mu(x)$ minimum.
    If $\bar\mu(x)\geq\delta/m$, then
    \[
        2=\sum_{y\in X}\bar\mu(y)\geq\abs{X} \cdot\bar\mu(x)\geq\abs{X}\cdot{\delta\over m}\quad\implies\quad\abs{X}\leq{2m\over\delta}.
    \]
    Otherwise, $\bar\mu(x)<\delta/m$.
    We then apply \Cref{lem:massbound} and use the inequalities $e^{-z/(1-z)}<1-z<e^{-z}$ for $0<z<1$ to bound
    \begin{align*}
        1 &\geq {1-{m\over\delta}\bar\mu(x)\over (1-\bar\mu(x))^{km}}>\exp\biggl\{{-{m\over\delta}\bar\mu(x)\over 1-{m\over\delta}\bar\mu(x)}+km\cdot\bar\mu(x)\biggr\}=\exp\biggl\{{m\cdot\bar\mu(x)\over 1-{m\over\delta}\bar\mu(x)}\biggl(k-{1\over\delta}-{km\over\delta}\bar\mu(x)\biggr)\biggr\},
    \end{align*}
    and so $\bar\mu(x)>{k\delta-1\over km}$.
    Therefore, since $k\delta\geq 2$,
    \[
        2=\sum_{y\in X}\bar\mu(y)>\abs{X} \cdot{k\delta-1\over km}\quad\implies\quad\abs{X} <{2km\over k\delta-1}.\qedhere
    \]
\end{proof}

\subsection{Cliques and even cycles}
In this section, we prove \Cref{thm:evencycle,thm:planarclique}.

We begin by computing $\optb(K_t,k)$.
\begin{theorem}\label{thm:optclique}
    For all $t\geq 2$ and all $k\geq 1$,
    \[
        \optb(K_t,k)={t\choose 2}^{-k{t\choose 2}}.
    \]
\end{theorem}
\begin{proof}
    The lower bound is realized by the uniform distribution on $E(K_t)$.
    \medskip

    For the upper bound, we have already shown that $\optb(K_2,k)=1$ (\Cref{prop:k2}), so we may suppose that $t\geq 3$.
    Fix any $\mu\in\Opt(K_t,k)$, which can be done thanks to \Cref{cor:achieved}.
    Note that $\abs{\supp\bar\mu}\geq t$ and that $\abs{\supp\mu}\geq{t\choose 2}$.

    Set $z=2/(t+1)$; we use a version of Bernoulli's inequality, $(1-x)^n<1-{nx\over 1+(n-1)x}$ for $0<x<1$ and $n>1$, to bound
    \begin{align*}
        (1-z)^{k{t\choose 2}} &\leq(1-z)^{{t\choose 2}}<1-{{t\choose 2}z\over 1+({t\choose 2}-1)z}=1-{{t\choose 2}\over t-1}z.
    \end{align*}
    Thus, thanks to \Cref{lem:massbound}, we know that $\bar\mu(x)>2/(t+1)$ for every $x\in\supp\bar\mu$.
    From here, we see that
    \[
        2=\sum_{x\in\supp\bar\mu}\bar\mu(x)>\abs{\supp\bar\mu}\cdot{2\over t+1}\quad\implies\quad\abs{\supp\bar\mu}<t+1\quad\implies\quad\abs{\supp\bar\mu}=t.
    \]
    Therefore, $\abs{\supp\mu}={t\choose 2}$, and so the claim follows from \Cref{prop:uniform}.
\end{proof}

Thus, the proof of \Cref{thm:planarclique} follows immediately from \Cref{lem:tumorb} (or \Cref{lem:tumorc} for $\bup{K_3}{1}$) and \Cref{thm:optclique}.
In fact, we have shown that
\[
    \numb_{\gcl C}(n,\bup{K_t}{k})={1\over(k!)^{{t\choose 2}}}\biggl({n\over{t\choose 2}}\biggr)^{k{t\choose 2}}+O(n^{k{t\choose 2}-k/(k+4)}),
\]
for all $t\geq 3$, $k\geq 1$ and $C\geq 2$.
\medskip

We next determine $\optb(C_4,k)$.
\begin{theorem}\label{thm:optc4}
    $\optb(C_4,k)=4^{-4k}$ for all $k\geq 1$.
\end{theorem}
\begin{proof}
    The lower bound is achieved by the uniform distribution on the edges of $C_4$.
    \medskip

    For the upper bound, fix any $\mu\in\Opt(C_4,k)$, which can be done thanks to \Cref{cor:achieved}.
    Set $X=\supp\bar\mu$; we claim that $\abs X=4$.
    Indeed, for any $x\in X$, \Cref{lem:massbound} tells us that
    \[
        1-2\bar\mu(x)\leq\bigl(1-\bar\mu(x)\bigr)^{4k}\leq\bigl(1-\bar\mu(x)\bigr)^4\quad\implies\quad\bar\mu(x)> 0.45.
    \]
    Therefore,
    \[
        2=\sum_{x\in X}\bar\mu(x)>0.45\cdot\abs X\quad\implies\quad \abs X<4.45,
    \]
    and so $\abs X=4$.
    We can therefore decompose ${X\choose 2}=\{e_1,f_1\}\cup\{e_2,f_2\}\cup\{e_3,f_3\}$ where $e_i,f_i$ are parallel edges, i.e.\ $e_i\cap f_i=\varnothing$.
    Since every copy of $C_4$ is uniquely determined by a pair of these parallel edges, we can write
    \begin{align*}
        \optb(\mu;C_4,k) &=\sum_{\{i,j\}\in{[3]\choose 2}}\mu(e_i)^k\mu(f_i)^k\mu(e_j)^k\mu(f_j)^k \leq \biggl(\sum_{\{i,j\}\in{[3]\choose 2}}\mu(e_i)\mu(f_i)\mu(e_j)\mu(f_j)\biggr)^k\\
                         &={1\over 2^k}\biggl(\biggl(\sum_{i=1}^3\mu(e_i)\mu(f_i)\biggr)^2-\sum_{i=1}^3\mu(e_i)^2\mu(f_i)^2\biggr)^k.
    \end{align*}
    We finally apply \Cref{cor:offdiag} and the AM--GM inequality to bound
    \begin{align*}
        \optb(\mu;C_4,k) &\leq {1\over 2^k}\biggl({1\over 8}\biggl(\sum_{i=1}^3\mu(e_i)\cdot\sum_{i=1}^3\mu(f_i)\biggr)^2\biggr)^k \leq{1\over 4^{2k}}\biggl({1\over 2}\sum_{i=1}^3\bigl(\mu(e_i)+\mu(f_i)\bigr)\biggr)^{4k}={1\over 4^{4k}}.\qedhere
    \end{align*}
\end{proof}

The proof of \Cref{thm:evencycle} now follows quickly.
\begin{proof}[Proof of \Cref{thm:evencycle}]
    The lower bounds are given in \cref{eqn:blowuplower}.

    Now, by applying \Cref{lem:tumorc}, we know that
    \[
        \numb_\plan(n,C_{2m})\leq\numb_{\gcl 3}(n,C_{2m})\leq\optb(C_m,1)\cdot n^m+O(n^{m-1/5}),
    \]
    for $m\geq 3$.
    Finally, \Cref{thm:optclique} gives $\optb(C_3,1)=3^{-3}$, \Cref{thm:optc4} gives $\optb(C_4,1)=4^{-4}$ and \Cref{thm:probb} gives $\optb(C_m,1)\leq 1/m!$ for all $m\geq 5$; hence the claim follows.
\end{proof}
\medskip

\subsection{Sufficiently large edge-blow-ups}\label{sec:largeblowup}
We conclude our study of $\optb(H,k)$ by proving \Cref{thm:blowupexact}.

\begin{theorem}\label{thm:largek}
    Let $H$ be a graph on $m$ edges with no isolated vertices and let $k$ be a positive integer.
    If $k\geq{\log(m+1)\over m\log(1+1/m)}$, then $\optb(H,k)=m^{-km}$.
\end{theorem}
\begin{proof}
    We begin by observing that if $k={\log(m+1)\over m\log(1+1/m)}$, then $(m+1)^{km-1}=m^{km}$.
    Since $k,m$ are positive integers and $m,m+1$ are coprime, this can happen only if $k=m=1$.
    This situation was covered in \Cref{prop:k2}, so we may suppose that $k>{\log(m+1)\over m\log(1+1/m)}$.

    Fix any $\mu\in\Opt(H,k)$, which can be done thanks to \Cref{cor:achieved} since $k>{\log(m+1)\over m\log(1+1/m)}\geq 1$.
    Set $x=1/(m+1)$ and observe that
    \begin{align*}
        (1-x)^{km} &< (1-x)^{{\log(m+1)\over\log(1+1/m)}}=\biggl({m\over m+1}\biggr)^{-{\log(m+1)\over\log(m/(m+1))}}={1\over m+1}=1-mx.
    \end{align*}
    Thus, thanks to \Cref{lem:massbound}, we see that $\mu(e)>1/(m+1)$ for every $e\in\supp\mu$.
    We conclude that
    \[
        1=\sum_{e\in\supp\mu}\mu(e)>{\abs{\supp\mu}\over m+1}\quad\implies\quad\abs{\supp\mu}<m+1\quad\implies\quad\abs{\supp\mu}=m,
    \]
    and so the claim follows from \Cref{prop:uniform}.
\end{proof}

The proof of \Cref{thm:blowupexact} then follows immediately from \Cref{lem:tumorb} and \Cref{thm:largek}.
\medskip

The lower bound of $k\geq{\log(m+1)\over m\log(1+1/m)}$ in \Cref{thm:largek} is tight for infinitely many graphs.
\begin{prop}\label{prop:edgetrans}
    Let $H$ be any edge-transitive graph on $m+1\geq 3$ edges.
    If $H^-$ is an $m$-edge subgraph of $H$ with no isolated vertices, then $\optb(H^-,k)>m^{-km}$ for all positive integers $k<{\log(m+1)\over m\log(1+1/m)}$.
\end{prop}
\begin{proof}
    First, note that ${\log(m+1)\over m\log(1+1/m)}>1$ since $m\geq 2$; hence the range for $k$ is nontrivial.

    Let $\mu$ denote the uniform distribution on $E(H)$.
    Since $H$ is edge-transitive, we know that $\numb(H,H^-)=m+1$ and so
    \begin{align*}
        {\optb(H^-,k)\over m^{-km}} &\geq{\optb(\mu;H^-,k)\over m^{-km}}=(m+1)\cdot\biggl({m\over m+1}\biggr)^{km}>(m+1)\cdot\biggl({m\over m+1}\biggr)^{-{\log(m+1)\over\log(m/(m+1))}}=1.\qedhere
    \end{align*}
\end{proof}

We remark that this is the reason that it is likely necessary to use the refined $\optb_\plan(H,k)$ mentioned in \Cref{rem:planarblowup} in order to determine $\numb_{\plan}(n,\bup{H}{k})$ for $H\in\{K_5^-,K_{3,3}^-\}$ and $k$ small.
For example, the proof of \Cref{prop:edgetrans} shows that $\optb(K_5^-,1)\geq 10^{-8}$, yet we think it is likely that $\optb_\plan(K_5^-,1)=9^{-9}$ since $K_5$ is not planar.

\section{Remarks and open problems}\label{sec:remarks}

The techniques introduced in this paper are far reaching.
Although we were able to compute $\optp(m)$ and $\optb(H,k)$ for certain $m$ and $H$, there is much we could not do.

\paragraph{Odd paths and even cycles}
The main question left open by this paper is that of determining $\optp(m)$ for $m\geq 4$.
\begin{conj}\label{conj:oddpath}
    For all $m\geq 2$, $\optp(m)$ is achieved by the uniform distribution on $E(C_m)$.
    In particular, $\optp(m)=8\cdot m^{-m}$.
\end{conj}
If true, then
\[
    \numb_\plan(n,P_{2m+1})=4m\biggl({n\over m}\biggr)^{m+1}+O(n^{m+4/5})\qquad\text{for all }m\geq 2,
\]
which would verify a conjecture of Ghosh et al.~\cite{ghosh_planarp5}, albeit with a worse error term than predicted.
Currently, we have only a proof for the cases of $m=2$ and $m=3$.

Even if \Cref{conj:oddpath} is true, the methods developed in this paper are likely too crude to achieve the posited error-term of $O(n^m)$, which would verify the conjecture of Ghosh et al.\ in full.
\medskip

Turning to even cycles, we conjecture the following:
\begin{conj}
    For all $m\geq 3$, $\optb(C_m,1)$ is achieved by the uniform distribution on $E(C_m)$.
    In particular, $\optb(C_m,1)=m^{-m}$.
\end{conj}
If true, then
\[
    \numb_\plan(n,C_{2m})=\biggl({n\over m}\biggr)^{m}+O(n^{m-1/5})\qquad\text{for all }m\geq 3.
\]
Currently, we have only a proof for the cases of $m=3$ and $m=4$.

It is likely that proving $\optb(C_m,1)=m^{-m}$ is well within reach for $m\in\{5,6\}$.
Indeed, for these values of $m$, one can use \Cref{lem:massbound} to show that $\optb(C_m,1)$ is achieved by a mass $\mu$ spanning exactly $m$ vertices.
Furthermore, one can show that $\bar\mu(x)=2/m$ for each $x\in\supp\bar\mu$.
We have not explored either of these cases any further.
Unfortunately, applying \Cref{lem:massbound} to $\optb(C_7,1)$ only allows us to say that this quantity is achieved by a mass spanning at most $8$ vertices.

\paragraph{Edge-blow-ups.}
The question of determining $\optb(H,k)$ is wide open for most graphs $H$.
One obvious lower-bound on $\optb(H,k)$ is the value achieved by the uniform distribution on $E(H)$.

\begin{question}
    For which graphs $H=(V,E)$ is $\optb(H,1)$ achieved by the uniform distribution on $E$?
    That is, for which graphs $H$ is $\optb(H,1)=\abs{E}^{-\abs E}$?
\end{question}
We have already noted that this is not the case for infinitely many graphs (\Cref{prop:edgetrans}).

Even though $\optb(H,1)$ is not always achieved by the uniform distribution on $E(H)$, it seems reasonable to expect that, given a finite set $X$, the quantity $\max\bigl\{\optb(\mu;H,1):\supp\mu\subseteq{X\choose 2}\bigr\}$ is achieved by the uniform distribution on the edges of \emph{some} graph.
If true, this leads to the following question, which could be interesting in its own right.
\begin{question}
    For a graph $H$ on $m$ edges with no isolated vertices, what bounds can be placed on the quantity
    \[
        \sup_G{\numb(G,H)\over\ \abs{E(G)}^m}?
    \]
\end{question}
Certainly this quantity is at least $1/m^m$ and is at most $1/m!$.
Additionally, we believe that the supremum can be replaced by a maximum unless $H=K_{1,m}$ or $H=mK_2$ for some $m\geq 2$ where $mK_2$ is the matching on $m$ edges.
\medskip

Finally, we still do not even know if $\optb(H,1)$ is achieved for many graphs.
Recall that $\optb(K_{1,m},1)$ and $\optb(mK_2,1)$ are never achieved for $m\geq 2$, yet $\optb(H,k)$ is achieved provided that $k\cdot\delta(H)\geq 2$ (\Cref{cor:achieved}).
\begin{conj}
    If $H$ is a graph with no isolated vertices and $k$ is a positive integer, then $\optb(H,k)$ is \emph{not} achieved if and only if $k=1$ and either $H=K_{1,m}$ or $H=mK_2$ for some $m\geq 2$.
\end{conj}

\paragraph{The reduction lemmas.}
Finally, we discuss the reduction lemmas in general.
First, as mentioned in \Cref{sec:reduce} after the statement of \Cref{lem:tumorc}, we believe the following to be true:
\begin{conj}
    Let $H$ be a graph on $m$ edges and let $k$ be a positive integer.
    If $k\cdot\delta(H)\geq 2$, then
    \[
        \numb_{\gcl C}(n,\bup{H}{k})={\optb(H,k)\over(k!)^m}\cdot n^{km}+o(n^{km}).
    \]
\end{conj}

Beyond this, it is natural to wonder if there is an analogous reduction lemma for even paths and odd cycles.
For example, the conjectured (asymptotic) extremal example for $\numb_\plan(n,P_{2m+2})$ is a modification of $\bup{C_m}{n/m}$ wherein a path is placed among the interior vertices of each blown-up edge (see~\cite[Conjecture 2]{ghosh_planarp5}); hence, we expect that the techniques used in this paper can be modified to tackle this question.
It is probably necessary to use more about the planar structure of the host-graph in order to extend the reduction lemmas to this situation.
\medskip

Interestingly, the reduction lemmas did not explicitly require the host-graph to have only linearly many edges.
By playing with the error-terms, one can extend each of the reduction lemmas to the collection of graphs $G$ which have no $K_{3,3}$ and $\abs{E(H)}\leq C\cdot\abs{V(H)}^{1+c}$ for each subgraph $H\subseteq G$, where $C>0$ is any fixed constant and $c>0$ depends on the particular situation at hand.
We opted to avoid this more general situation for the sake of readability.

Furthermore, it was not crucial that the host-graph avoided copies of $K_{3,3}$.
Indeed each of the reduction lemmas can be reworked to handle the case when the host-graph avoids copies of $K_{3,t}$ for some fixed $t\geq 3$.
In particular, the reduction lemmas apply to the class graphs which can be embedded onto any surface of a fixed genus.
However, the fact that one side of this forbidden biclique has size $3$ appears to be necessary for each of our arguments.
It seems unlikely that similar reduction lemmas could be pushed through if the host-graph only avoids copies of, say, $K_{4,4}$.

\medskip

Finally, it is pertinent to point out that the techniques developed in this paper can likely be extended to prove stability results for $\numb_\plan(n,H)$ for various graphs $H$.
This would, however, likely require a few new ideas.

\bibliographystyle{abbrv}
\bibliography{references}

\end{document}